\documentclass{elsarticle}

\usepackage{setspace}
\usepackage{amsmath}
\usepackage{amsfonts}
\usepackage{mathrsfs}
\usepackage{bm}
\usepackage[title]{appendix}
\usepackage{graphicx}
\usepackage{tikz}
\usepackage{pgfplots}
\usepackage{caption}
\usepackage{subcaption}
\pgfplotsset{compat = 1.3}
\usepackage[margin=1in]{geometry}

\journal{Journal of Computational Physics}

\begin{document}

\begin{frontmatter}

\title{High-order accurate finite difference discretisations on fully unstructured dual quadrilateral meshes}

\author[rvt2,rvt3]{Y.~Pan\fnref{fn1}\corref{cor1}}
\ead{yllpan@berkeley.edu}

\author[rvt2,rvt3]{P.-O.~Persson\fnref{fn3}}
\ead{persson@berkeley.edu}

\address[rvt2]{Department of Mathematics, University of California, Berkeley, Berkeley, CA 94720, United States}
\address[rvt3]{Mathematics Group, Lawrence Berkeley National Laboratory, 1 Cyclotron Road, Berkeley, CA 94720, United States}
\cortext[cor1]{Corresponding author}
\fntext[fn1]{Graduate student, Department of Mathematics, University of California, Berkeley}
\fntext[fn3]{Professor, Department of Mathematics, University of California, Berkeley}

\begin{keyword}
 Finite Differences, %
 High-order methods, %
 Unstructured meshes
\end{keyword}

\begin{abstract}
    We present a novel approach for high-order accurate numerical differentiation on unstructured meshes of quadrilateral elements. To differentiate a given function, an auxiliary function with greater smoothness properties is defined which when differentiated provides the derivatives of the original function. The method generalises traditional finite difference methods to meshes of arbitrary topology in any number of dimensions for any order of derivative and accuracy. We demonstrate the accuracy of the numerical scheme using dual quadrilateral meshes and a refinement method based on subdivision surfaces. The scheme is applied to the solution of a range of partial differential equations, including both linear and nonlinear, second and fourth order equations, and a time-dependent first order equation.
\end{abstract}

\end{frontmatter}
	
\section{Introduction}

Some of the most efficient methods for numerical solution of partial differential equations are based on finite difference (FD) techniques. These include the standard textbook methods for elliptic, parabolic, and hyperbolic equations \cite{leveque_fd} as well as more sophisticated numerical schemes such as the compact FD method \cite{lele_compact,visbal02compactfd} or WENO method \cite{shu99wenotriangular,shu_siamreview}. The schemes can be highly computationally efficient, due to the point-wise stencils and sparse connectivities, as well as being highly robust in the context of under-resolved features such as shocks or other discontinuities. However, most of these methods are only applicable to structured grids, and for real-world problems that require unstructured meshes and boundary-fitted elements, finite elements and finite volume methods are more commonly used.

Many approaches have been proposed for extending FD methods to unstructured grids. The multi-block methods \cite{carpenter_99,mattson_10,nordstrom09multiblock,rai_86,zhang04multiblock} are based on applying standard finite difference schemes inside (large) blocks of structured nodes, and strategies such as upwinding, interpolation and numerical fluxes used to connect the blocks. A related class of schemes are the overset methods \cite{chesshire90,jespersen1997overflow,steger1983chimera} popular in CFD applications, which interpolate between structured blocks without requiring conforming interfaces. Finite difference methods have been derived for adaptive structured grids, in particular using octree grids in the so-called AMR methods \cite{bell_amr,colella_amr}. Finally, high-order unstructured finite difference/volume stencils can been derived using higher-dimensional Taylor expansions \cite{barth1990higher,jensen72,liszka96,liszka80,ollivier2002volume,perrone75}. While many of these methods can be competitive, it remains fair to say that the challenges that arise from unstructured meshes remain the main reason why finite difference methods are not as widely used in practical applications.

In this work, we describe a new way to formulate finite difference stencils on arbitrary unstructured quadrilateral meshes. We first present an approach for deriving 1D stencils on grids of irregular nodes. Our strategy is to introduce a carefully chosen auxiliary function with greater smoothness properties than the original function. This function can then be differentiated using standard techniques, from which derivatives of the original function may be extracted. While these resulting derivative stencils could also be obtained by standard Taylor expansion techniques, the introduction of a new function with higher regularity enable the application of a range of traditional finite difference techniques on regular grids.

To extend the method to 2D and higher dimensions, we refine the original quadrilateral mesh using a subdivision-based scheme, and apply our 1D stencils on the (unstructured) grid lines from the dual mesh. This isolates the irregular nodes, around which we form specialized stencils by excluding the nodes in one of the quadrants. We demonstrate the resulting finite difference stencils by verifying the orders of convergence on both structured and unstructured meshes under refinement. Finally, we apply the techniques on several partial differential equations, including Poisson's equation, the biharmonic equation, and the (nonlinear) minimal surface equation. We also demonstrate the method on a time-dependent problem of advection equation discretised using upwinding.

\section{1D Numerical Differentiation}

\subsection{Basics}
A brief summary on basic numerical differentiation is provided. For a function $f(x) \in C^{\infty}(a,b) \subset \mathbb{R}$, we can expand the function as a Taylor series centred at some point $x_0 \in (a,b)$
\begin{equation}
    f(x_0 + h) = f(x_0) + hf'(x_0) + \sum_{k=2}^\infty \frac{h^k}{k!}f^{(k)}(x_0)
\end{equation}
In practice instead of the function we have an ordered sequence of points $\{ x_i \in [a,b],~ i=1,2,...,n ,~ x_i < x_{i+1},~ x_1 = a,~ x_n = b \}$ known as a grid, on which the function is evaluated. In the case where the spacing between all points $x_{i+1} - x_i = h$, the grid is said to be regular, otherwise the grid is irregular.

On a regular grid using Taylor expansions as written above we can approximate the first derivative of $f(x_i)$ using any of the following

\begin{align}
    f'(x_i) &= \frac{1}{h} \big( f(x_{i+1}) - f(x_i) \big) + O(h) \label{eq:1d_forward},~~ i < n \\
    f'(x_i) &= \frac{1}{h} \big( f(x_i) - f(x_{i-1}) \big) + O(h) \label{eq:1d_backward},~~ i > 1 \\
    f'(x_i) &= \frac{1}{2h} \big( f(x_{i+1}) - f(x_{i-1} ) \big) + O(h^2) \label{eq:1d_centre},~~ 1 < i < n
\end{align}
Equations \ref{eq:1d_forward} and \ref{eq:1d_backward} are known as forward and backward differences respectively whilst Equation \ref{eq:1d_centre} is known as the centred difference. The forward and backward differences in Equations \ref{eq:1d_forward} and \ref{eq:1d_backward} are first order accurate, whilst the central difference formula in Equation \ref{eq:1d_centre} is second order accurate. 

By taking into accounts more points on the grid, these can be extended to compute higher order accurate approximations to the first derivative or to compute higher derivatives. For instance, the well known 2nd order central difference approximation of the 2nd derivative is given as
\begin{equation}
    \label{eq:1d_centre2}
    f''(x_i) = \frac{1}{h^2}\big( f(x_{i+1}) -2f(x_i) + f(x_{i-1}) \big) + O(h^2)
\end{equation}

We can also rewrite each of the above Equations \ref{eq:1d_forward} - \ref{eq:1d_centre2} in stencil form. For instance, Equation \ref{eq:1d_centre} can be rewritten as
\begin{equation}
    f'(x_i) = \frac{1}{h}
    \begin{bmatrix}
        -1/2 \\
        1/2
    \end{bmatrix} \cdot 
    \begin{bmatrix}
        f(x_{i-1}) \\
        f(x_{i+1})
    \end{bmatrix} + O(h^2), ~~ 1<i<n
\end{equation}
where the left vector is known as the weights/coefficients. The two vectors together form what is known as the stencil of the approximation.

\subsection{Irregular grids}
\label{sect:irregular1d}
Consider a function $f(x) \in C^{\infty}(a,b) \subset \mathbb{R}$, but this time the function is instead to be evaluated on an irregular grid $\{ x_1,...,x_n \}$. Without loss of generality, we may consider the example of finding $f'(x)$ at a point $x_k \in (a,b)$ where $x_{k} - x_{k-1} = h_1,~ x_{k+1} - x_{k} = h_2$ for $h_1, h_2 \in \mathbb{R}^+$ positive real numbers.

To numerically differentiate the function at $x_k$, we can parametrise $x$ by a new variable $\xi$ using an interpolatory spline of order at least $m \geq \frac{n}{2}$. The parametrisation may be chosen such that the grid $\{ \xi_i, i=1,2,...,n,~ \xi_1 = 0,~ \xi_n = 1 \}$ in $\xi$-space is regular and evenly spaced; that the spline function is interpolatory implies that $x(\xi_i) = x_i$. Composing the function $f(x)$ with the parametrisation $x(\xi)$ we get
\begin{equation}
    \tilde{f}(\xi) = f\big(x (\xi) \big)
\end{equation}
where in general the function $\tilde{f}(\xi) \in C^\infty \big( (\xi_1,\xi_n) \big)$. Taking derivatives of this function $\tilde{f}(\xi)$
\begin{equation}
    \tilde{f}'(\xi) = f'(x) x'(\xi)
\end{equation}
we see that given the derivative of $\tilde{f}(\xi)$, the derivative $f'(x)$ can be easily computed, as the Jacobian $x'(\xi)$ can be approximated using Equations \ref{eq:1d_forward}-\ref{eq:1d_centre}.

Computing $\tilde{f}'(\xi)$ is however in general more challenging than simply applying Equations \ref{eq:1d_forward}-\ref{eq:1d_centre}. Defining $h = \max\{ x_{k}-x_{k-1},~ x_{k+1}-x_{k} \}$, applying the centred difference in Equation \ref{eq:1d_centre} to $\tilde{f}(\xi)$ at $\xi_k$, we find that even as $h$ decreases under refinement of the grid, that the error of the approximation remains constant. To see why this is the case, we can write Equation \ref{eq:1d_centre} as
\begin{equation}
    \tilde{f}'(\xi_k) = \frac{1}{2h} \big( \tilde{f}(\xi_{k+1}) - \tilde{f}(\xi_{k-1} ) \big) + \frac{h^2}{6}\tilde{f}'''(\xi_k) + O(h^3)
\end{equation}
Computing $\tilde{f}'''(\xi)$ at $\xi_k$ to analyse the leading error term
\begin{align*}
    \tilde{f}'(\xi_k) &= f'(x_k) x'(\xi_k) \\
    \tilde{f}''(\xi_k) &= f''(x_k) x'(\xi_k)^2 + f'(x_k) x''(\xi_k) \\
    \tilde{f}'''(\xi_k) &= f'''(x_k) x'(\xi_k)^3 + 3f''(x_k)x'(\xi_k)x''(\xi_k) + f'(x_k) x'''(\xi_k)
\end{align*}
we note that $x''(\xi_k) = O(\frac{1}{h}), x'''(\xi_k) = O(\frac{1}{h^2}), x^{(l)}(\xi_k) = O(\frac{1}{h^{l-1}}), l \leq m$. This behaviour can be explained by the following argument: as $x(\xi)$ is chosen to be a spline, the curve is refined via insertion of new knots which does not alter the value of $x(\xi_k)$ or its first derivative $x'(\xi_k)$ as $h$ decreases. However by definition of the derivative as $h$ decreases the second derivative $x''(\xi_k)$ must then increase proportionally. Repeating this argument for higher derivatives establishes the claim above. Thus for the leading error term we have
\begin{equation}
    \tilde{f}'(\xi_k) = \frac{1}{2h} \big( \tilde{f}(\xi_{k+1}) - \tilde{f}(\xi_{k-1} ) \big) + \frac{h^2}{6} \underbrace{{\tilde{f}'''(\xi_k)}}_{O(\frac{1}{h^2})}
\end{equation}
which remains constant with decreasing $h$. This is unless $\tilde{f}'(\xi_k) = 0$, in which case would imply
\begin{align*}
    \tilde{f}'(\xi_k) &= 0 \\
    \tilde{f}'''(\xi_k) &= f'''(x_k) x'(\xi_k)^3 + 3f''(x_k)x'(\xi_k)x''(\xi_k) = O(\frac{1}{h}) \\
    \tilde{f}'(\xi_k) &= \frac{1}{2h} \big( \tilde{f}(\xi_{k+1}) - \tilde{f}(\xi_{k-1} ) \big) + \frac{h^2}{6} \underbrace{{\tilde{f}'''(\xi_k)}}_{O(\frac{1}{h})}
\end{align*}
and first order accuracy approximating $\tilde{f}'(\xi_k)$ using Equation \ref{eq:1d_centre}. By the same reasoning, if $\tilde{f}(\xi)$ were such that $\tilde{f}'(\xi_k) = 0, \tilde{f}''(\xi_k) = 0$ the method would become second order accurate.

Observing this we look to construct a function $\tilde{F}(\xi)$ such that $\tilde{F}^{(l)}(\xi) = 0, l=1,2,...$ and from which we can infer information on the derivatives $f^{(l)}(x_k)$. Specifically we define the set of functions 
\begin{align}
    F_l(x; x_k) &= f( x ) - f( x_k ) - \sum_{j=1}^l \frac{( x - x_k )^j}{j!} f^{(j)}(x_k) \\
    \tilde{F}_l(\xi; x_k) &= f( x(\xi) ) - f( x_k ) - \sum_{j=1}^l \frac{( x(\xi) - x_k )^j}{j!} f^{(j)}(x_k)
\end{align}
where $l = 1,2,..$. The function is constructed by subtracting the terms of the Taylor series of $f(x)$ centred at $x_k$ up to order $l$, such that for a given value of $l$, we have $F^{(j)}(x_k; x_k) = \tilde{F}^{(j)}(\xi_k; x_k) = 0$ for $j=0,1,2,...,l$.

The trick now is to use the knowledge that the derivatives of $\tilde{F}_l(\xi_k;x_k)$ are zero when applying difference formulae to approximate the derivatives of $\tilde{F}_l(\xi_k;x_k)$. In general, this gives us a linear system which we can solve for the values of $f^{(m)}(x_k), m = 1,2,...,l$. To see how this works, we can apply the central difference formula in Equation \ref{eq:1d_centre} to approximate the derivative of the function $\tilde{F}_1(\xi;x_k)$ at $\xi_k$
\begin{align*}
    0 = F_1^{'}(x_k;x_k) = \tilde{F}_1^{'}(\xi_k;x_k) = \frac{1}{2h} \big( \tilde{F}_k(\xi_{k+1}&) - \tilde{F}_k(\xi_{k-1}) \big) \\
    0 = \tilde{f}( \xi_{k+1} ) - \big( x(\xi_{k+1})-x_k \big)f'(x_k) - \tilde{f}( \xi_{k-1} ) + \big( &x(\xi_{k-1})-x_k \big)f'(x_k) \\
    f'(x_k) = \frac{\tilde{f}( \xi_{k+1} ) - \tilde{f}( \xi_{k-1} )}{x(\xi_{k+1}) - x(\xi_{k-1})}
\end{align*}
This is of course the well known first order accurate central difference approximation for irregular spaced grids in 1D, which is consistent with the order of accuracy suggested in the discussion above.

\subsection{2nd order accurate approximation of 1st derivative}
To build intuition for the method, we consider the example of finding a second order accurate approximation of $f'(x)$ at $x=x_k$ with the parametrisation $x(\xi) \in C^m (\xi_1,\xi_n), m \geq \frac{n}{2}$. To do this we consider the function $\tilde{F}_2(\xi;x_k)$

\begin{equation}
    \tilde{F}_2(\xi;x_k) = F_2\big( x(\xi);x_k \big) = f\big(x(\xi)\big) - f(x_k) - (x(\xi) - x_k)f'(x_k) - \frac{(x(\xi) - x_k)^2f''(x_k)}{2}
\end{equation}
It can easily be checked that $F_2'(x_k;x_k) = F_2''(x_k;x_k) = \tilde{F}_2'(\xi_k;x_k) = \tilde{F}_2''(\xi_k;x_k) = 0$. Combining this with Equations \ref{eq:1d_centre} and \ref{eq:1d_centre2} to approximate the first and second derivatives of $\tilde{F}(\xi;\xi_k)$ at $\xi = \xi_k$, we obtain the system of equations
\begin{equation*}
    \frac{1}{2} \tilde{F}_2(\xi_{k+1};x_k) - \frac{1}{2} \tilde{F}_2(x_{k-1};x_k) = 0, \quad \tilde{F}_2(\xi_{k+1}; x_k) -2\tilde{F}_2(\xi_k;x_k) + \tilde{F}_2(\xi_{k-1};x_k) = 0 \\
\end{equation*}
This time we have two equations and two unknowns to be solved for: $f'(x_k), f''(x_k)$. Denoting $x_{k-1} - x_k = h_l, x_{k+1}-x_k = h_r$, the linear system can be written as

\begin{equation}
    \begin{split}
        \begin{bmatrix}
            -\frac{1}{2} & \frac{1}{2}\\
            1 & 1
        \end{bmatrix}
        \begin{bmatrix}
            h_l & \frac{ h_l^2 }{2}\\
            h_r & \frac{ h_r^2 }{2}
        \end{bmatrix}
        \begin{bmatrix}
            f'(x_k) \\
            f''(x_k)
        \end{bmatrix}
            =
        \begin{bmatrix}
            -\frac{1}{2} & 0 & \frac{1}{2} \\
            1 & -2 & 1
        \end{bmatrix}
        \begin{bmatrix}
            f(x_{k-1}) \\
            f(x_k) \\
            f(x_{k+1})
        \end{bmatrix} \\
        \begin{bmatrix}
            \frac{-h_l + h_r}{2} & -\frac{ h_l^2 }{4} + \frac{ h_r^2 }{4}\\
            h_l + h_r & \frac{ h_1^2 }{2} + \frac{ h_r^2 }{2}
        \end{bmatrix}
        \begin{bmatrix}
            f'(x_k) \\
            f''(x_k)
        \end{bmatrix}
        =
        \begin{bmatrix}
            -\frac{1}{2} & 0 & \frac{1}{2} \\
            1 & -2 & 1
        \end{bmatrix}
        \begin{bmatrix}
            f(x_{k-1}) \\
            f(x_k) \\
            f(x_{k+1})
        \end{bmatrix} \\
        \begin{bmatrix}
            f'(x_k) \\
            f''(x_k)
        \end{bmatrix}
        =
        \frac{1}{h_r - h_l}
        \begin{bmatrix}
            \frac{h_r}{h_l} & \frac{h_l^2 - h_r^2 }{h_l h_r} & -\frac{h_l}{h_r}\\
            -\frac{2}{h_l} & \frac{2h_r - 2h_l}{h_l h_r} & \frac{2}{h_r}
        \end{bmatrix}
        \begin{bmatrix}
            f(x_{k-1}) \\
            f(x_k) \\
            f(x_{k+1})
        \end{bmatrix}
    \end{split}
\end{equation}
Solving this system we get both an approximation to the second derivative $f''(x_k)$ in addition to an approximation for the first derivative $f'(x_k)$. Specifically for a general $h_l, h_r$ the above expressions are the well known three point stencils for a second order approximation to the first derivative and a first order approximation to the second derivative on an irregular grid. In the special case of a regular grid where $ x_{k}-x_{k-1} = x_{k+1} - x_{k} = h $, the system simplifies to
\begin{equation}
    \begin{split}
        \begin{bmatrix}
            -\frac{1}{2} & \frac{1}{2} \\
            1 & 1
        \end{bmatrix}
        \begin{bmatrix}
            -h & \frac{h^2}{2} \\
            h & \frac{h^2}{2}
        \end{bmatrix}
        \begin{bmatrix}
            f'(x_k) \\
            f''(x_k)
        \end{bmatrix}
        =
        \begin{bmatrix}
            -\frac{1}{2} & 0 & \frac{1}{2} \\
            1 & -2 & 1
        \end{bmatrix}
        \begin{bmatrix}
            f(x_{k-1}) \\
            f(x_k) \\
            f(x_{k+1})
        \end{bmatrix} \\
        \begin{bmatrix}
            h & 0 \\
            0 & h^2
        \end{bmatrix}
        \begin{bmatrix}
            f'(x_k) \\
            f''(x_k)
        \end{bmatrix}
        =
        \begin{bmatrix}
            -\frac{1}{2} & 0 & \frac{1}{2} \\
            1 & -2 & 1
        \end{bmatrix}
        \begin{bmatrix}
            f(x_{k-1}) \\
            f(x_k) \\
            f(x_{k+1})
        \end{bmatrix} \\
        \begin{bmatrix}
            f'(x_k) \\
            f''(x_k)
        \end{bmatrix}
        =
        \begin{bmatrix}
            -\frac{1}{2h} & 0 & \frac{1}{2h} \\
            \frac{1}{h^2} & -\frac{2}{h^2} & \frac{1}{h^2}
        \end{bmatrix}
        \begin{bmatrix}
            f(x_{k-1}) \\
            f(x_k) \\
            f(x_{k+1})
        \end{bmatrix}
    \end{split}
\end{equation}
and we recover the well known three point second order accurate regular grid centred difference stencils for the first and second derivatives.

\subsection{rth order accurate approximation of qth derivative}
We now consider the general case of generating a rth order accurate approximation to the qth derivative of $f(x)$ on an irregular grid $\{ x_1,...,x_n \}$ at $x_k$ parametrised by $x=x(\xi) \in C^m (\xi_1,\xi_n), m \geq \frac{n}{2}$. To do this we look at the function $F_p(x;x_k)$ where $p = r+q-1$

\begin{equation}
    F_p(x; x_k) = f( x ) - f( x_k ) - \sum_{j=1}^p \frac{( x - x_k )^j}{j!} f^{(j)}(x_k)
\end{equation}

We have that $F^{(l)}_p(x_k; x_k) = 0$ for $l=0,1,2,...,p$. As a result we may now write $p$ difference equations for approximating each of the $p$ derivatives of $\tilde{F}_p(\xi; x_k) = F_p\big( x(\xi); x_k \big)$ at $\xi = \xi_k$, with the knowledge that the derivatives are all equal to zero.

For the $p$ difference equations, there is some freedom as to which stencil can be used. The stencil used for the qth derivative of $\tilde{F}_p(\xi; x_k)$ should at least be rth order accurate. For the other equations, the rule to be followed is that for each order derivative that is lower than $q$, the stencil used must be an order higher, and vice versa. For instance, the stencil for the equation for the $q-3$ derivative of $\tilde{F}_p(\xi;x_k)$ must be at least $r+3$ order accurate, and the stencil for the equation for the $q+2$ derivative of $\tilde{F}_p(\xi;x_k)$ must be at least $r-2$ order accurate.

In the prior subsection, the stencils chosen were 2nd order accurate for both 1st and 2nd derivatives, which satisfy the above condition. The stencil chosen for the 2nd derivative could in fact be changed to be only 1st order accurate, and would not affect the rate of convergence of either of the approximations to the derivatives $f'(x_k), f''(x_k)$.

We may now solve for the $p$ unknowns $f'(x_k), f''(x_k), ..., f^{(p)}(x_k)$. Denoting the stencil points $\{ x_1, ..., x_k,...,x_n \}$ and their corresponding stencil weights for the lth derivative as $c^{(l)}_j$, this amounts to solving a linear system of the form

\begin{equation} 
    \label{eq:1d_linear_system}
    CXDu = \bar{C}f
\end{equation}

where the terms are as follows:
\begin{enumerate}
    \item $f$ is a vector with $n$ entries, $f_{j} = f(x_j)$,
    \item $\bar{C}$ is a $p \times n$ matrix, $C_{ij} = c^{(i)}_j$, matrix of coefficient weights,
    \item $C$ is a $p \times (n-1)$ matrix, which is formed by deleting the $k$th column of $\bar{C}$,
    \item $X$ is a $(n-1) \times p$ Vandermonde matrix, $X_{ij} = ( x_i - x_k )^j $, of the positions of the stencils points
    \item $D$ is a $p \times p$ diagonal matrix, $D_{ij} = \delta_{ij} \cdot \frac{1}{j!}$, a constant matrix of the factorial terms in the denominator of the Taylor series of $f$
    \item $u$ is a vector with $p$ entries, $u_{i} = f^{(i)}(x_k)$, containing the desired derivatives of $f$
\end{enumerate}
The right hand side to the equation is the result of directly applying a regular grid stencil with step size $h=1$ on an irregular grid, the result of which is corrected by the left hand side mesh correction term $CXD$. As seen in the example of the second order accurate approximation to the first derivative in the previous section, in the case where the grid is in fact regular, the mesh correction becomes diagonal. Specifically, the jth diagonal entry of $CXD$ for a regular grid is equal to $\frac{1}{h^j}$ recovering the standard finite difference stencils.

\subsection{Local truncation error}
\label{sect:lte}
To analyse the local truncation error from the above procedure on an irregular grid, we consider the case of an approximation of the qth derivative of $\tilde{f}(\xi)$ with order r accuracy denoted $\tilde{f}^{(q)}_{\xi,r}$ obtained using a difference stencil. As the function $\tilde{f}(\xi)$ is defined in $\xi$-space discretised using a regular grid, the error in the approximation equals
\begin{equation}
    \tilde{f}^{(q)}(\xi) - \tilde{f}^{(q)}_{\xi,r} = Dh^d\tilde{f}^{(q+r)}(\xi) + O(h^{r+1})
\end{equation}
where $D$ is a constant. For the leading error term, we expand the derivative of $\tilde{f}(\xi)$ as
\begin{equation}
    \begin{split}
        \tilde{f}^{(m)}(\xi) &= \sum_{j=1}^m D_j f^{(j)}\big( x(\xi) \big) \cdot G_{m,j}(\xi)\\
        G_{m,j}(\xi) &= O(1/h^{(m-j)})
    \end{split}
\end{equation}
where $D_j$ are some constants. This equality is established using numerical induction. In the case $m=1$, the claim holds as $x'(\xi) = O(1)$ from the discussion above in Sect. \ref{sect:irregular1d}. For the inductive step we write
\begin{align*}
    \tilde{f}^{(m+1)}(\xi) &= \frac{d}{d\xi} \sum_{j=1}^m D_j f^{(j)} \big( x(\xi) \big) \cdot G_{m,j}(\xi) \\
    &= \sum_{j=1}^m D_j \bigg[ f^{(j+1)}\big( x(\xi) \big) x'(\xi) \cdot G_{m,j}( \xi ) + f^{(j)}\big( x(\xi) \big) \cdot G_{m,j}'(\xi) \bigg]
\end{align*}
Where for the last term
\begin{align*}
    G_{m,j}'(\xi) = G_{m+1,j}(\xi) &= O \bigg( \frac{G_{m,j}(\xi+h) - G_{m,j}(\xi)}{h} \bigg) \\
    &= \frac{O(1/h^{m-j})}{h} = O(1/h^{m+1-j})
\end{align*}
which proves the claim. Using this we can examine the order of the leading error term to the approximation
\begin{equation}
    \begin{split}
        h^df^{(q+r)}(\xi) &= \sum_{j=p+1}^{q+r} D_j f^{(j)}\big( x(\xi) \big) \cdot G_{q+r,j}(\xi) \\
        &= O( h^r / h^{q+r-p-1} ) = O( h^{p+1-q} )
    \end{split}
\end{equation}
where $p$ is the largest integer such that $f^{(1)}(x(\xi)) = f^{(2)}(x(\xi)) = ... = f^{(p)}(x(\xi)) = 0$. This motivates the choice of $\tilde{F}_p(\xi;x)$ in the previous section for a rth order accurate approximation of the qth derivative where $p = r+q-1$.

\subsection{Preconditioning the Vandermonde system}
Solving Equation \ref{eq:1d_linear_system} requires the solution of a Vandermonde type system where in practice the Vandermonde matrix can have more columns than rows. The shape of this matrix is due to the fact that in general there are more points used for approximating the derivatives than the number of derivatives approximated. It is however well-known that Vandermonde type systems can be very ill-conditioned and thus solving the system can become difficult to perform accurately.

To solve this system a right preconditioner can be applied such that
\begin{equation}
    CX\bigg(K^{-1}K\bigg)Du = \bar{C}f
\end{equation}
where $K^{-1}_{ij} = \delta_{ij}\frac{1}{x_j^j}$ is chosen to be the diagonal Jacobi preconditioner acting on the columns of $X$. The rationale behind this choice lies in the structure of the Vandermonde matrix, wherein the jth column of the matrix is of the form $( x_1^j, ~...,~x_i^j )^T$ and so for large $j$ the values of the column can become extremely large or small depending on the sizes of $x_i$. In our case the entries $x_i$ denote distances between nodes of the mesh and so are of similar magnitude. Multiplying the column by $\frac{1}{x_j^j}$ thus scales the column such that the values are of magnitude $O(1)$ and allows for much better conditioning of the Vandermonde system.

\subsection{Numerical differentiation example}
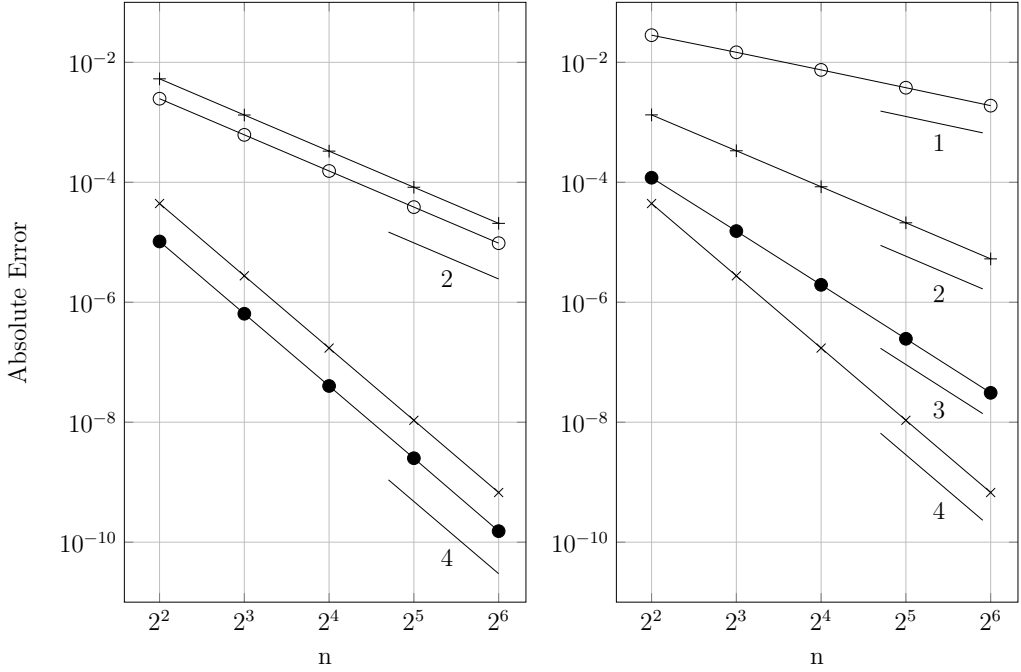
\begin{figure}
    \centering
    \begin{tikzpicture}[font=\large, scale=0.8]
        \begin{loglogaxis}[
            width=0.5\textwidth,
            height=0.7\textwidth,
            xlabel = {n},
            xlabel style={ yshift=-1ex },
            ylabel = {Absolute Error},
            ylabel style={ yshift=2ex },
            ymin=1e-11, ymax=1e-1,
            xmin=3, xmax=80,
            xtick={4,8,16,32,64},
            xticklabels={$2^2$,$2^3$,$2^4$,$2^5$,$2^6$},
            ytick={1e-2, 1e-4, 1e-6, 1e-8, 1e-10},
            yticklabels={$10^{-2}$,$10^{-4}$,$10^{-6}$,$10^{-8}$,$10^{-10}$},
            grid = both,
            grid style = {line width=.1pt, draw=gray!15},
            major grid style = {line width=.2pt, draw=gray!50},
        ]
        \addplot[mark=x, color=black, mark size=3pt]
        table[ x=n, y=reg1 ]{dat/diff1d_reg.csv};
        
        \addplot[mark=*, color=black, mark size=3pt]
        table[ x=n, y=reg2 ]{dat/diff1d_reg.csv};
        
        \addplot[domain=26:64, samples=2] {0.01*x^-2};
        \node at (axis cs:42, 2.5e-6) {2};
        
        \addplot[mark=+, color=black, mark size=3pt]
        table[ x=n, y=reg3 ]{dat/diff1d_reg.csv};
        
        \addplot[mark=o, color=black, mark size=3pt]
        table[ x=n, y=reg4 ]{dat/diff1d_reg.csv};
        
        \addplot[domain=26:64, samples=2] {0.0005*x^-4};
        \node at (axis cs:42, 5.5e-11) {4};
        
        \end{loglogaxis}
    \end{tikzpicture}%
    ~
    \begin{tikzpicture}[font=\large, scale=0.8]
        \begin{loglogaxis}[
            width=0.5\textwidth,
            height=0.7\textwidth,
            xlabel = {n},
            xlabel style={ yshift=-1ex },
            ymin=1e-11, ymax=1e-1,
            xmin=3, xmax=80,
            xtick={4,8,16,32,64},
            xticklabels={$2^2$,$2^3$,$2^4$,$2^5$,$2^6$},
            ytick={1e-2, 1e-4, 1e-6, 1e-8, 1e-10},
            yticklabels={$10^{-2}$,$10^{-4}$,$10^{-6}$,$10^{-8}$,$10^{-10}$},
            grid = both,
            grid style = {line width=.1pt, draw=gray!15},
            major grid style = {line width=.2pt, draw=gray!50},
        ]
        \addplot[mark=x, color=black, mark size=3pt]
        table[ x=n, y=reg1 ]{dat/diff1d_reg.csv};
        \label{plot:diff1d_d1}
        
        \addplot[domain=26:60, samples=2] {0.04*x^-1};
        \node at (axis cs:42, 4.5e-4) {1};
        
        \addplot[mark=*, color=black, mark size=3pt]
        table[ x=n, y=irreg2 ]{dat/diff1d_irreg.csv};
        \label{plot:diff1d_d2}
        
        \addplot[domain=26:60, samples=2] {0.006*x^-2};
        \node at (axis cs:42, 1.4e-6) {2};
        
        \addplot[mark=+, color=black, mark size=3pt]
        table[ x=n, y=irreg3 ]{dat/diff1d_irreg.csv};
        \label{plot:diff1d_d3}
        
        \addplot[domain=26:60, samples=2] {0.003*x^-3};
        \node at (axis cs:42, 1.6e-8) {3};
        
        \addplot[mark=o, color=black, mark size=3pt]
        table[ x=n, y=irreg4 ]{dat/diff1d_irreg.csv};
        \label{plot:diff1d_d4}
        
        \addplot[domain=26:60, samples=2] {0.003*x^-4};
        \node at (axis cs:42, 3.3e-10) {4};
        
        \end{loglogaxis}
    \end{tikzpicture}
    \caption{Numerical differentiation of $f(x)=e^{0.7x} - x^2$ at $x=0$ on regular grid (left) and irregular grid (right). \ref{plot:diff1d_d1} denotes the absolute error in the 1st derivative, \ref{plot:diff1d_d2} denotes the absolute error in the 2nd derivative, \ref{plot:diff1d_d3} denotes the absolute error in the 3rd derivative, \ref{plot:diff1d_d4} denotes the absolute error in the 4th derivative. }
    \label{fig:diff1d}
\end{figure}

As an example of we compute a fourth order accurate approximation of the first derivative of the function $f(x) = e^{0.7x} + x^2$ at the point $x=0$. Following the discussion in the previous subsection, we look at the function
\begin{equation}
    F_4(x;0) = f(x) - f(0) - \sum_{j=1}^4 \frac{x^j}{j!}f^{(j)}(0)
\end{equation}
We consider the regular grid $P_{\text{reg}} = \{ -2/n, -1/n, 0, 1/n, 2/n \}$ and the irregular grid $P_{\text{irreg}}= \{ -2/n, -1/n, 0, 2/5n, 4/5n  \}$ where $n = 1,2,...$ on which to perform numerical differentiation. For both $P_{\text{reg}}, P_{\text{irreg}}$ we have the matrices
\begin{align*}
    \bar{C} =  \begin{bmatrix}
    \frac{1}{12} & \frac{-2}{3} & 0 & \frac{2}{3} & \frac{-1}{12}\\
    \frac{-1}{12} & \frac{4}{3} & \frac{-5}{2} & \frac{4}{3} & \frac{-1}{12}\\
    \frac{-1}{2} & 1 & 0 & -1 & \frac{1}{2}\\
    1 & -4 & 6 & -4 & 1
    \end{bmatrix}
    ,
    C =  \begin{bmatrix}
    \frac{1}{12} & -\frac{2}{3} & \frac{2}{3} & -\frac{1}{12}\\
    -\frac{1}{12} & \frac{4}{3} & \frac{4}{3} & -\frac{1}{12}\\
    -\frac{1}{2} & 1 & -1 & \frac{1}{2}\\
    1 & -4 & -4 & 1
    \end{bmatrix}
\end{align*}
where each row in $\bar{C}$ is simply the well known 5-point centred difference stencils for the first to fourth derivatives with step size $h=1$, and $C$ is formed by deleting the third column of $\bar{C}$ corresponding to the point at $x=0$, and
\begin{align*}
    D =  \begin{bmatrix}
    \frac{1}{1!} & 0 & 0 & 0\\
    0 & \frac{1}{2!} & 0 & 0\\
    0 & 0 & \frac{1}{3!} & 0\\
    0 & 0 & 0 & \frac{1}{4!}
    \end{bmatrix}
\end{align*}
the diagonal matrix corresponding to the factorial denominators in the Taylor series expansion. The difference in the two cases lie in the $X$ term in Eq. \ref{eq:1d_linear_system}
\begin{align*}
    X_{\text{reg}} =  \begin{bmatrix}
    \frac{-2}{n} & \frac{4}{n^2} & \frac{-8}{n^3} & \frac{16}{n^4}\\
    \frac{-1}{n} & \frac{1}{n^2} & \frac{-1}{n^3} & \frac{1}{n^4}\\
    \frac{1}{n} & \frac{1}{n^2} & \frac{1}{n^3} & \frac{1}{n^4}\\
    \frac{2}{n} & \frac{4}{n^2} & \frac{8}{n^3} & \frac{16}{n^4}
    \end{bmatrix} ,
    X_{\text{irreg}} =  \begin{bmatrix}
    \frac{-2}{n} & \frac{4}{n^2} & \frac{-8}{n^3} & \frac{16}{n^4}\\
    \frac{-1}{n} & \frac{1}{n^2} & \frac{-1}{n^3} & \frac{1}{n^4}\\
    \frac{2}{5n} & \frac{4}{25n^2} & \frac{8}{125n^3} & \frac{16}{625n^4}\\
    \frac{4}{5n} & \frac{16}{25n^2} & \frac{64}{125n^3} & \frac{256}{625n^4}
    \end{bmatrix}
\end{align*}
and the vectors $f$
\begin{align*}
    f_{\text{reg}} = 
    \begin{bmatrix}
    f(\frac{-2}{n}) \\
    f(\frac{-1}{n}) \\
    f(0) \\
    f(\frac{1}{n}) \\
    f(\frac{2}{n})
    \end{bmatrix} ,
    f_{\text{irreg}} =
    \begin{bmatrix}
    f(\frac{-1}{n}) \\
    f(\frac{-1}{2n}) \\
    f(0) \\
    f(\frac{2}{5n}) \\
    f(\frac{4}{5n})
    \end{bmatrix}
\end{align*}
It can be verified that for the regular grid $CX_{\text{reg}}D$ is a diagonal matrix where the $j$-th diagonal entry is $\frac{1}{n^j}$. Solving the linear system $CXDu = \bar{C}f$ for an approximation to the derivatives of $f(x) = e^{0.7x} + x^2$ at $x=0$ using the two grids, we obtain the results shown in Fig. \ref{fig:diff1d}. For the regular grid, the errors in 1st and 2nd derivatives converge with fourth order accuracy, whilst the errors in 3rd and 4th derivatives converge with 2nd order accuracy. For the irregular grid we observe that the error in the 1st derivative is 4th order accurate, and that for each higher derivative the order of accuracy decreases by one as expected. Furthermore whilst the accuracy of approximations on the irregular grid to the 2nd and 4th derivatives are lower compared to those on the regular grid as a result of the decreased order of convergence, the accuracy of the approximations to the 1st and 3rd derivatives is very close to that on the regular grid.

\section{Numerical Differentiation in Higher Dimensions}
\subsection{Primal vs. dual meshing}
To differentiate a function $f$ in $d$-dimensions, a grid of points is commonly specified as $d$-dimensional array of points on which the function is evaluated. For instance in 2D, a regular grid can be specified as $\{ x_{i,j} ~|~ i<n, j<m \}$, where $x_{i+1,j} - x_{i+1,j} = h_1,~ x_{i,j+1} - x_{i,j} = h_2$. However in dimensions greater than one, in addition to variable spacing between grid points there is the added complexity of the topology on which a function may be defined. While the formulation here is general for any number of dimensions, we focus here only on the cases of 2- and 3- dimensions as we are mainly interested in applications in these cases.

For many applications, instead of a structured grid, an unstructured mesh $\mathscr{M}$ is instead given, loosely defined as a collection of points in $d$-dimensional space known as vertices and of connections between two vertices. A closed set of edges is known as a face. The number of edges incident to a vertex is known as the valency of the vertex.

We focus here on unstructured quadrilateral meshes, wherein most faces on a mesh have exactly four sides. The mesh of a regular grid is a special case of a quad mesh wherein every interior vertex has a valency of four. In practice meshes in 2-dimensions may exhibit topological defects in the form of interior vertices with a valency unequal four, known as extraordinary points, or faces consisting of number of edges unequal to four, known as an extraordinary face.

\begin{figure}
    \centering
    \includegraphics[scale=0.45]{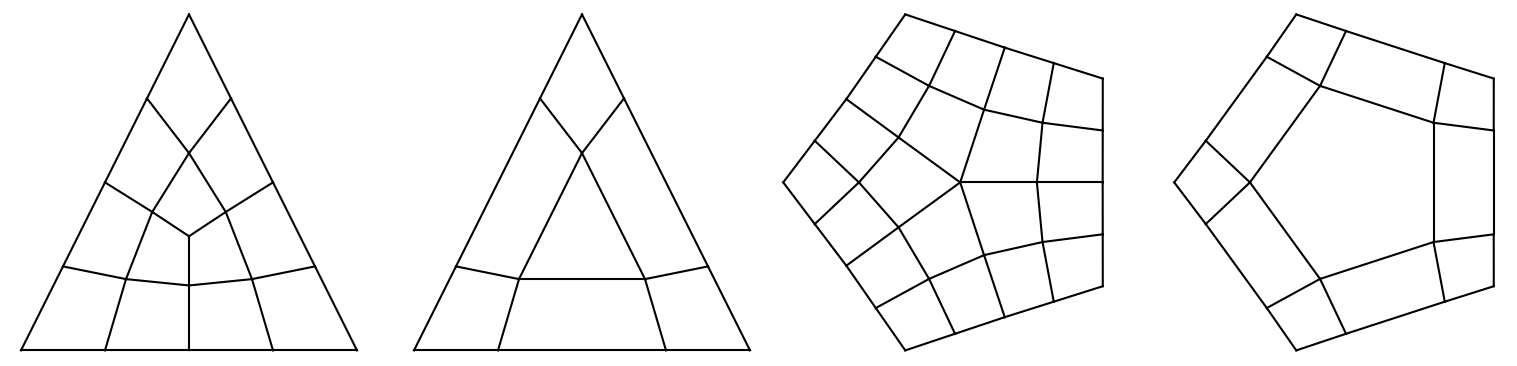}
    \caption{Comparison of primal quad meshes and their duals. The leftmost mesh has a central extraordinary vertex with a valency of three, the edges incident to the vertex are deleted to form its dual shown in centre left. Primal mesh in centre right has a central extraordinary vertex with a valency of five, the edges incident to this vertex are deleted to form its dual on the right.}
    \label{fig:primal_dual}
\end{figure}

In the case where all mesh topological defects are represented in the form of extraordinary points, the mesh is known as a primal quad mesh, whilst a mesh where all topological defects are represented as extraordinary faces is known as a dual quad mesh. The two types of representation are in fact equivalent and it is easy to convert a primal quad mesh to its dual form and vice versa, as shown in Fig. \ref{fig:primal_dual}.

3-dimensional quad meshes can be similarly classified, where a primal 3-dimensional quad mesh is one where topological defects are represented using extraordinary points with a valency unequal to six. Likewise a 3-dimensional dual quad mesh is one where topological defects are represented using extraordinary faces where the number of edges is unequal to four. As with the 2-dimensional case a similar isomorphism exists between 3-dimensional primal and dual forms of a mesh.

In this text all meshes considered are quad meshes in their dual forms. All meshes were generated first in their primal forms before being converted into dual forms with the procedure shown in Fig. \ref{fig:primal_dual}. Laplacian smoothing was then applied with a convergence threshold of $h \cdot 10^{-4}$, with $h$ denoting the average initial edge length of the mesh.

\subsection{Dual mesh refinement}
One notable example of dual mesh refinement in 2-dimensions is Doo-Sabin subdivision \cite{doo_sabin}, based on knot insertion on uniform biquadratic B-Spline surfaces. Classical Doo-Sabin refinement however does not preserve boundaries under refinement. To remedy this we modify Doo-Sabin subdivision by following the procedure used by Catmull and Clark \cite{catmull_clark}. While we focus here only on the refinement process for two-dimensional dual surface meshes, we note that it can be extended to higher dimensions by extending the procedure outlined by Catmull and Clark \cite{catmull_clark} as performed here.

Standard subdivision schemes such as Doo-Sabin or Catmull-Clark are derived from knot insertion on uniform biquadratic and bicubic B-Spline surfaces respectively. B-Splines and their tensor product surfaces are parametrised using knot vectors, where subdivision refinement is equivalent to inserting a new knot at the midpoint of each interval between successive knots. For our modification, our refinement procedure is equivalent to inserting two knots in every interval between successive knots of a biquadratic B-Spline surface, at one-third and two-thirds of the interval splitting the interval into three equal parts.

To begin the refinement step, we first compute for every face on the mesh its midpoint by averaging the position of all its vertices. For every edge, its midpoint is calculated by averaging the positions of its two endpoints. These midpoints do not form part of the refined mesh, but are to be used to calculate new points of the refined mesh.

\begin{figure}
    \centering
    \includegraphics[scale=0.35]{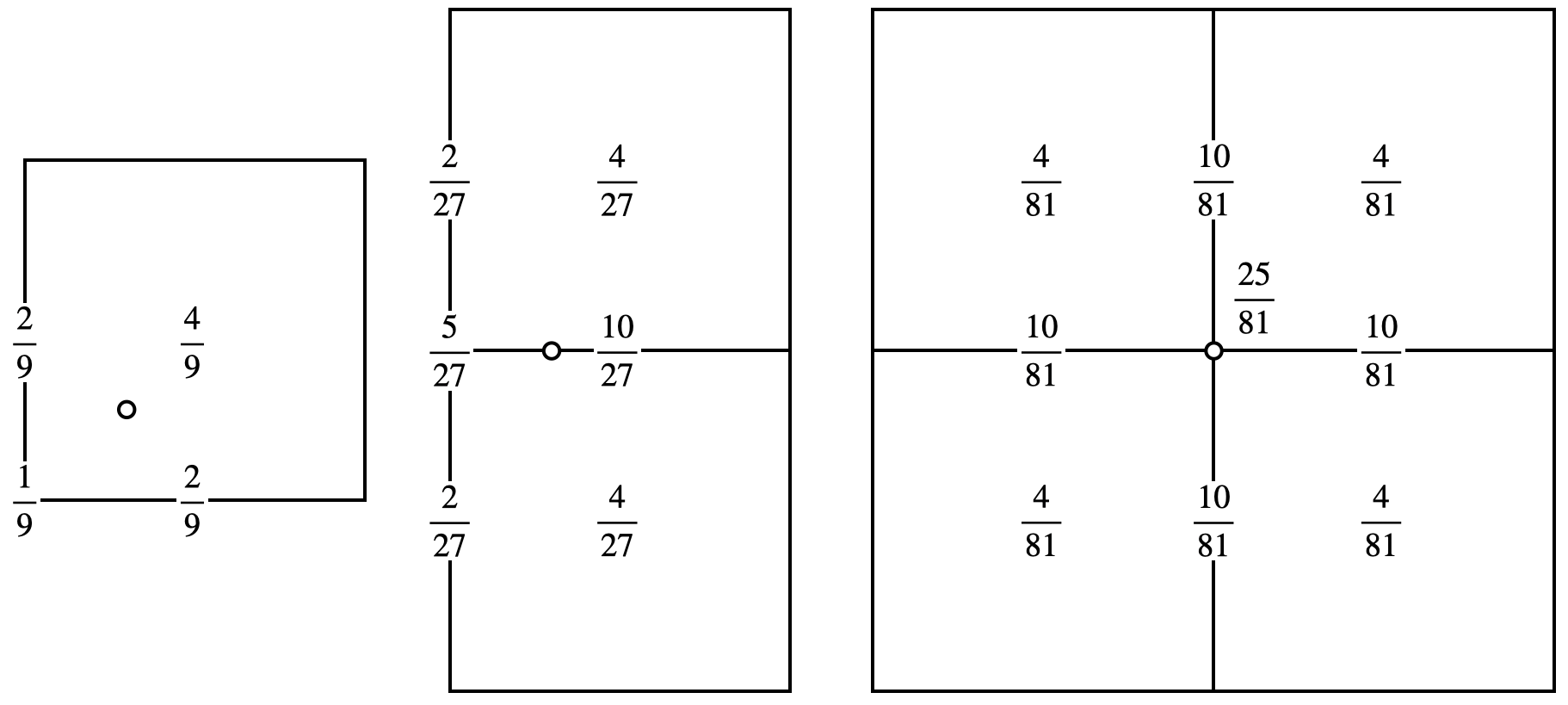}
    \caption{Weights for dual mesh refinement. From left to right: face points, edge points, vertex points.}
    \label{fig:refinement_rules}
\end{figure}

\begin{figure}
    \centering
    \includegraphics[scale=0.4]{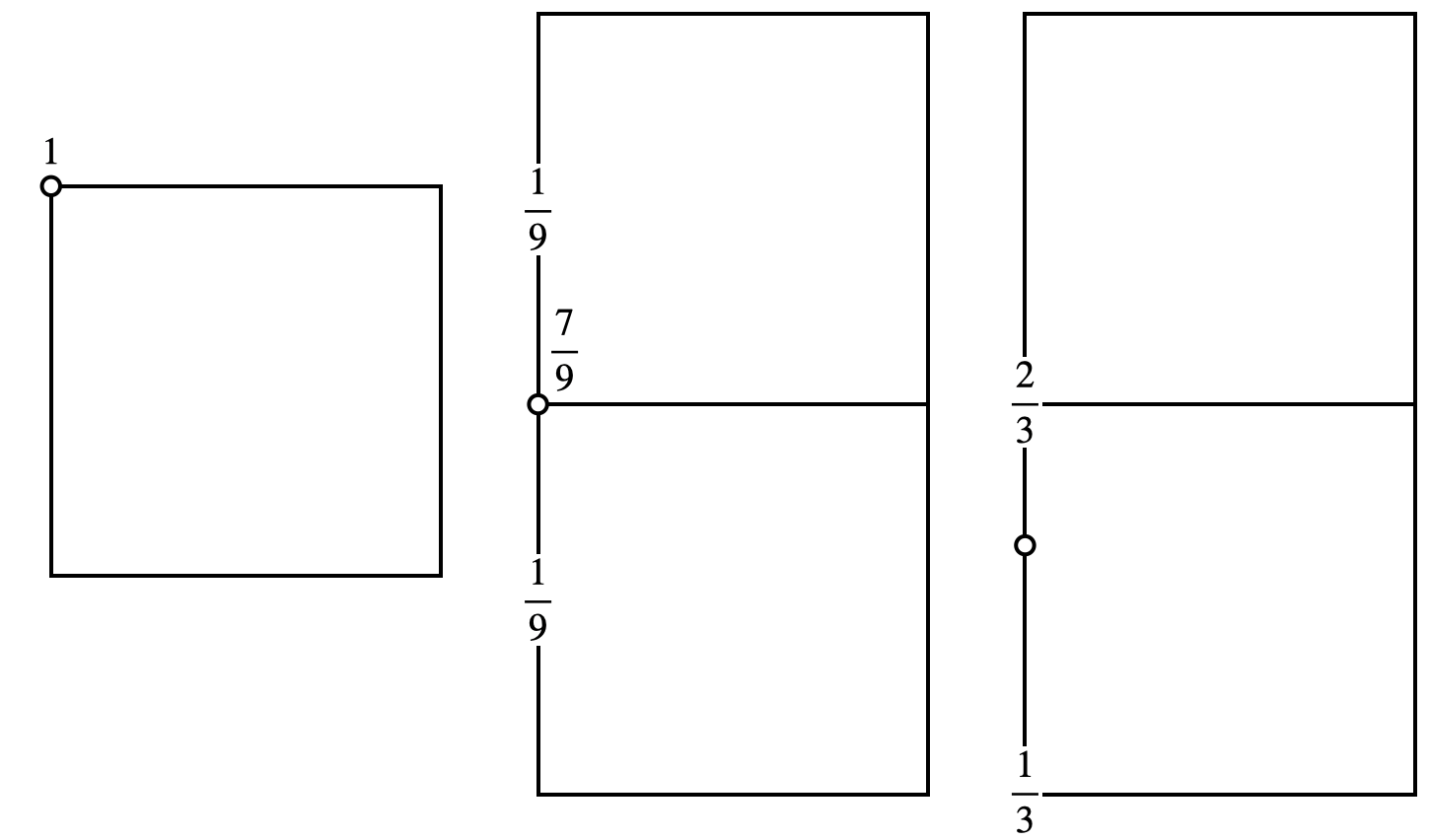}
    \caption{Modified weights for refinement on boundaries. From left to right: vertex point for corners, vertex points on boundaries, edge points on boundaries.}
    \label{fig:refinement_rules_boundary}
\end{figure}

Using the above midpoints, new points for the refined mesh are calculated by taking weighted sums of the midpoints and vertices on the mesh. All weights are shown in Fig. \ref{fig:refinement_rules}.
\begin{enumerate}
    \item For every face, one new vertex is created per vertex defining the face. These vertices are known as face points and are calculated by taking a weighted sum of the vertex, the face midpoint, and the two midpoints of edges incident to the vertex adjacent to the face.
    \item For every edge, one new vertex is created for each endpoint of the edge. These vertices are known as edge points and are calculated by taking a weighted sum of the endpoint, midpoints of the two faces adjacent to the edge, and midpoints of all edges incident to the endpoint adjacent to those two faces.
    \item For every vertex, a new vertex is created, known as a vertex point. The vertex point is a weighted sum of the original vertex, and midpoints of all edges incident to the vertex and midpoints of all faces adjacent to the vertex.
\end{enumerate}
To ensure that boundaries are preserved, a modification is made to the weights of edge and vertex points on boundaries, shown in Fig. \ref{fig:refinement_rules_boundary}. Corner vertex points are set to the position of the original vertex, while boundary vertex points are set to be an weighted average the original vertex and only the midpoints of the two incident boundary edges. The last modification is to boundary edge points, which are modified to be a weighted average of its two endpoints only.

\begin{figure}
    \centering
    \includegraphics[scale=0.56]{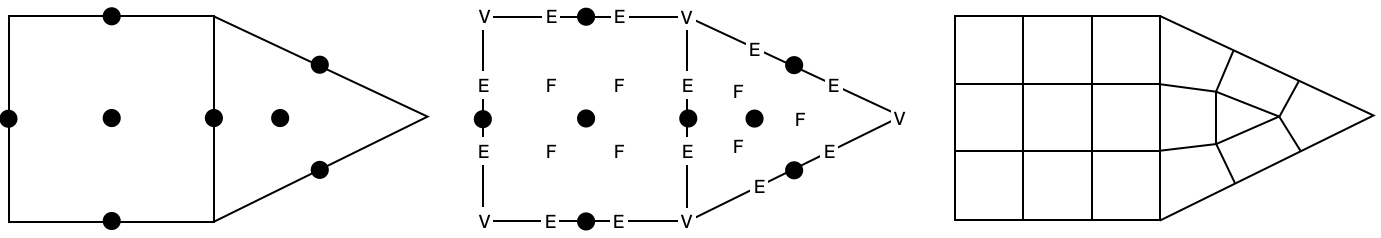}
    \caption{Dual mesh refinement procedure. From left to right: midpoints of every face and edge are computed shown in bold, new vertex/edge/face points are computed using the precomputed midpoints and stencils shown in Figs. \ref{fig:refinement_rules}-\ref{fig:refinement_rules_boundary}, new edges are drawn between newly generated vertex/edge/face points to complete refinement step.}
    \label{fig:refinement_procedure}
\end{figure}

\begin{figure}
    \centering
    \includegraphics[scale=0.45]{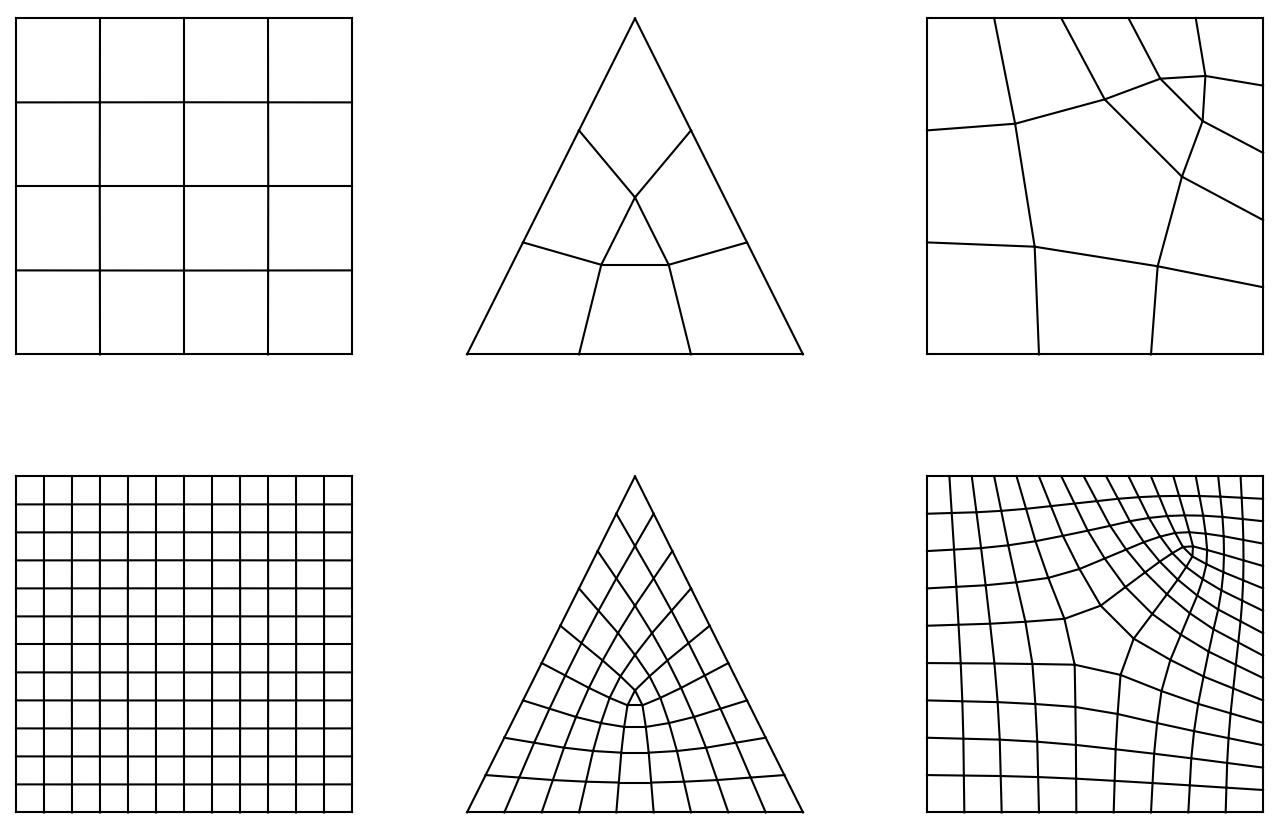}
    \caption{Dual mesh refinement examples.}
    \label{fig:refinement_examples}
\end{figure}

To complete the refinement step, edges are drawn between each new vertex points and its four new adjacent edge points, each edge point and its two adjacent face point, each edge point and its adjacent edge point, and each face point to its two adjacent face points. This is shown in Fig. \ref{fig:refinement_procedure} and some refinement examples are shown in Fig. \ref{fig:refinement_examples}. In contrast to many standard refinement procedures where each refinement roughly halves the edge lengths, edge lengths under this procedure decrease roughly by a factor of three.

\subsection{Numerical differentiation on dual meshes}
One key property of dual quad meshes is that for every interior vertex, the valence is equal to exactly $2d$ where $d$ denotes the dimension. This allows us for every vertex to locally identify $d$ lines that intersect exactly at the vertex. Identification of these lines is performed heuristically: edges are paired off such the lines they form always intersect at the vertex. Extending this procedure to every vertex in the mesh allows for the definition of curves $d$ curves  $\{ c_1(\xi_1), c_2(\xi_2), ..., c_d(\xi_d) \}$ at every vertex which provide a local parametrisation for the computation of derivatives.

\begin{figure}
    \centering
    \includegraphics[scale=0.45]{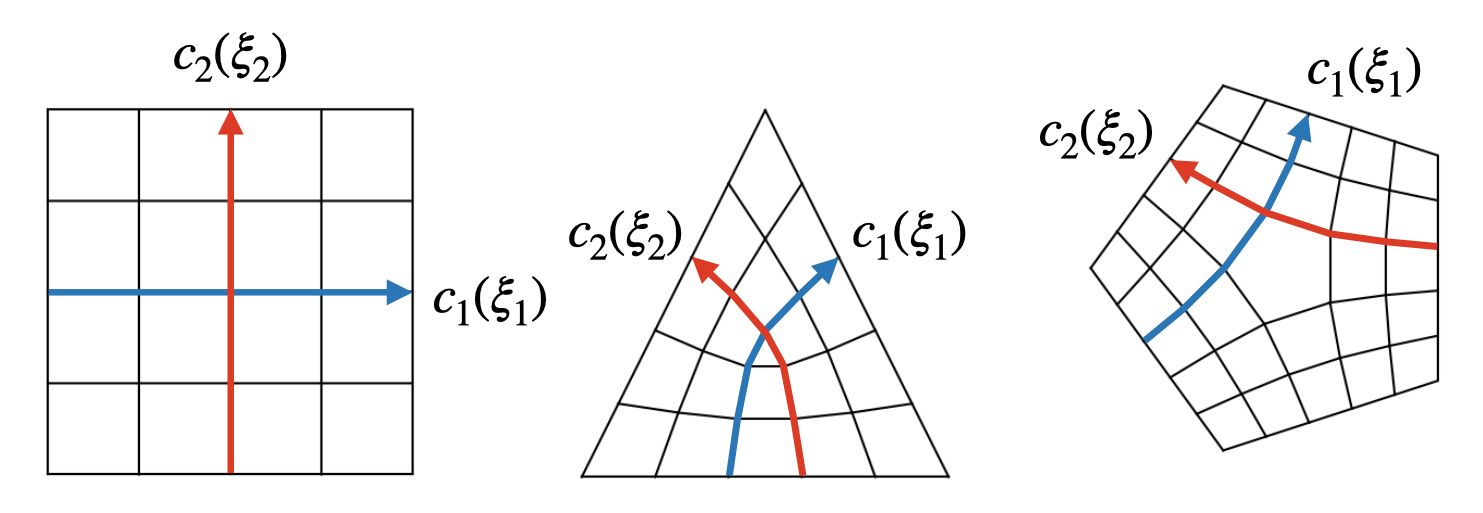}
    \caption{Curve identification for 2D dual quad meshes.}
    \label{fig:lines_2d}
\end{figure}

To differentiate a function $f(x_1,...,x_d)$ at a vertex $(x_{1k},...,x_{dk})$ of the mesh, we define an auxiliary function $\tilde{F_p}$ on the local parametrisation  $(\xi_1,...,\xi_d)$ such that
\begin{enumerate}
    \item Derivatives of $f$ can be approximated via numerical differentiation of $\tilde{F_p}$. This is easier to do than directly differentiating $f$ as $\tilde{F_p}$ is a function on the parametric space $(\xi_1,...,\xi_d)$ which is discretised using a regular uniform grid.
    \item $\frac{\partial^j \tilde{F}_p}{\partial \xi_1^{j_1} ... \partial \xi_d^{j_d}} = 0,$ where $\sum_{i=1}^d j_i = j $ for $j = 0,1,...,p$
\end{enumerate}
Analogously with the one dimensional case, $\tilde{F_p}$ is constructed via subtracting from $f$ terms of the Taylor series of $f$ centred at $(x_{1k},...,x_{dk})$ with total derivative order less than or equal to $p$. In two dimensions, relabelling $x_1, x_2$ as $x,y$ the function becomes
\begin{equation}
    \begin{split}
        F_p(x,y;x_{k},y_{k}) = f(x,y) - \sum_{j=0}^p &\sum_{l=0}^{p-j} \frac{( x - x_{k} )^j}{j!} \frac{( y - y_{k} )^l}{l!} \frac{\partial^p f}{\partial x^j \partial y^l} (x_{k},y_{k}) \\
        \tilde{F}_p(\xi_1, \xi_2; x_{k},y_{k}) &= F_p\big( x(\xi_1,\xi_2), x(\xi_1,\xi_2) \big)
    \end{split}
\end{equation}
and in three dimensions relabelling $x_1, x_2, x_3$ as $x,y,z$
\begin{equation}
    \begin{split}
        F_p(x&,y,z;x_{k},y_{k},z_{k}) = f(x,y,z) \\
        - \sum_{j=0}^p &\sum_{l=0}^{p-j} \sum_{m=0}^{p-j-l} \frac{( x - x_{k} )^j}{j!} \frac{( y - y_{k} )^l}{l!} \frac{( z - z_{k} )^m}{m!} \frac{\partial^p f}{\partial x^j \partial y^l \partial z^m} (x_{k},y_{k},z_{k}) \\
        \tilde{F}_p(\xi_1, \xi_2 &, \xi_3 ; x_{k},y_{k},z_{k}) = F_p\big( x(\xi_1, \xi_2, \xi_3), y(\xi_1, \xi_2, \xi_3), z(\xi_1, \xi_2, \xi_3) \big)
    \end{split}
\end{equation}
To compute a rth order accurate of the derivative $\frac{\partial^q f}{\partial x_1^{q_1} ... \partial x_d^{q_d} }$ where $\sum_{i=1}^d q_i = q$, it must be chosen that $p = r+q-1$. This choice is justified by extension of the argument outlined in Sect. \ref{sect:lte}.

Finally we apply regular grid numerical differentiation stencils for each of the derivatives $\frac{\partial^j}{\partial \xi_1^{j_1} ... \partial \xi_d^{j_d}}, j \leq p$ to $\tilde{F_p}$ at the vertex $\{ x_{1k},...,x_{dk} \}$, with the knowledge that each of the derivatives of $\tilde{F}_p$ at the vertex are equal to zero. The rule for stencil selection remains that the chosen stencil for the derivative $\frac{\partial^j}{\partial \xi_1^{j_1} ... \partial \xi_d^{j_d}}$ must be at least $r-j+1$ order accurate. This sets up a linear system which we can solve to obtain the derivatives $\frac{\partial^j}{\partial x_1^{j_1} ... \partial x_d^{j_d}}$ of $f$ at $\{ x_{1k},...,x_{dk} \}$. This again results in a linear system of the form
\begin{equation} \label{eq:stencil_sys}
    CXDu = \bar{C}f
\end{equation}
mirroring the the one dimensional case where $C$ is a matrix of the stencils weights omitting the column corresponding to the point $(x_k,y_k,z_k)$ at which the derivative is to be computed, $X$ a $d$-dimensional Vandermonde matrix of the positions of stencil points, $D$ a diagonal matrix of the factorial denominators from the Taylor Series terms, and $\bar{C}f$ the result of applying regular grid stencils with step size $h=1$ directly to the values of the function $f$ at the stencil points. In the case where the mesh is in fact a regular grid, the left hand side $CXD$ becomes diagonal recovering the standard regular grid difference stencils.

\subsection{Stencils at topological defects}
\begin{figure}
    \centering
    \includegraphics[scale=0.4]{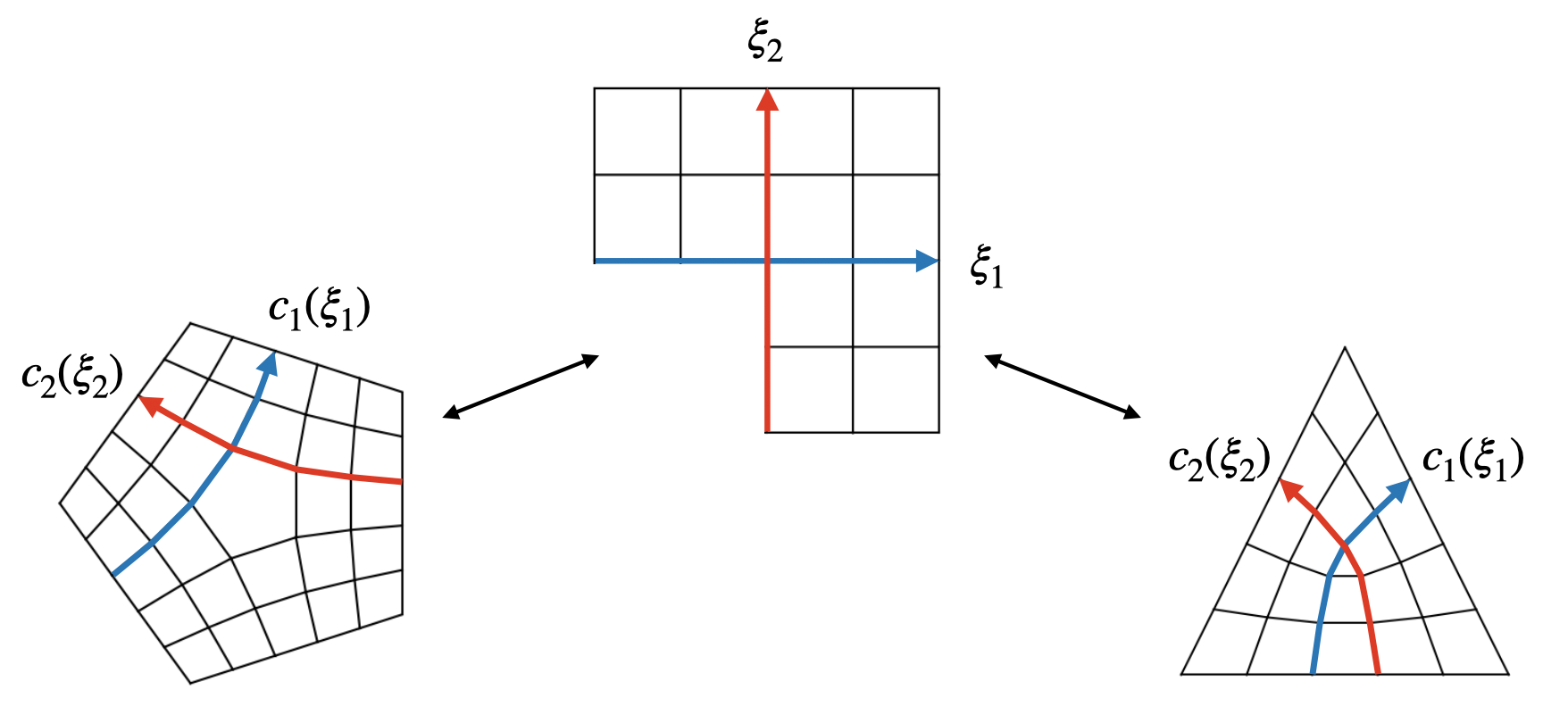}
    \caption{Local three quadrant parametrisation around extraordinary face.}
    \label{fig:param_3quad}
\end{figure}

An extra challenge is presented at vertices close to topological defects, in that even with the definition of $d$ curves through the vertex there does not exist a simple mapping locally of the region to a Cartesian grid. This means that points for the regular grid stencil cannot be uniquely chosen as a result of extraordinary faces in the mesh. This is circumvented by omitting one portion of the parametric space at the extraordinary face such that for the remaining space there exists in a simple mapping to $2^d - 1$ quadrants (in 2D)/octants (in 3D) Cartesian space. For two dimensions this is shown in Fig. \ref{fig:param_3quad}.

\begin{figure}
    \centering
    \includegraphics[scale=0.45]{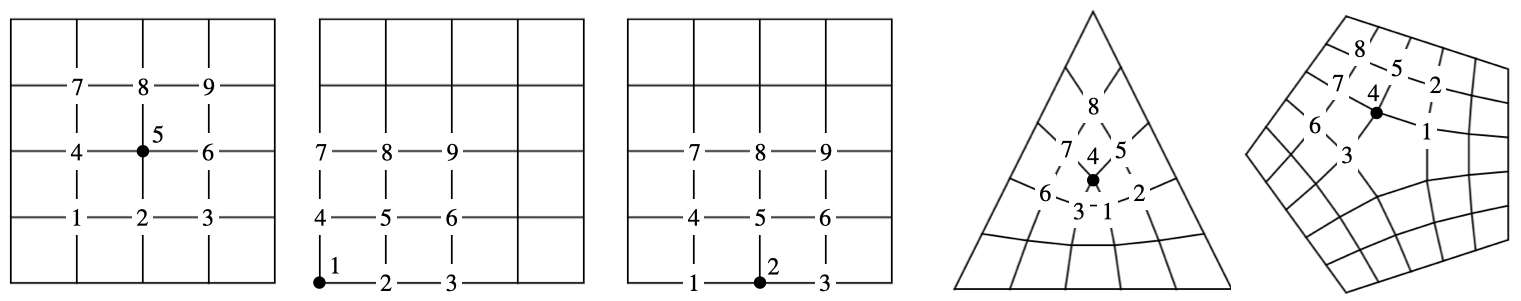}
    \caption{Stencil points and their numbering for second order accurate central difference $\frac{\partial f}{\partial x}, \frac{\partial f}{\partial y}$ example. Solid coloured vertex indicates point $(x_k,y_k)$ at which the function derivative is approximated. From left to right: symmetric stencil for regular point, one sided stencil for boundary point, one sided stencil for corner point, symmetric 3 quadrant stencil at an extraordinary face on a triangular mesh, symmetric 3 quadrant stencil at an extraordinary face on a pentagonal mesh.}
    \label{fig:stencil_2d_centre3}
\end{figure}

Stencils omitting one portion of parametric space may then be chosen for each of the $\frac{\partial^j}{\partial \xi_1^{j_1} ... \partial \xi_d^{j_d}}$ derivatives of $\tilde{F}_p$ such that they are $k-j+1$ order accurate. This requirement means that extraordinary faces on the mesh cannot be too close to one another, with the minimum distance set by the desired order of accuracy and the type of stencil chosen. As an example, for a second order accurate symmetric stencil for $\frac{\partial f}{\partial x}$, eights stencil points are needed to guarantee second order accuracy, and so extraordinary faces must be separated by at least one ordinary face. The points chosen in this case are shown in Fig. \ref{fig:stencil_2d_centre3}.

\subsection{Efficient implementation of stencil assembly}
The stencil assembly procedure described above requires the assembly of the matrix system Equation \ref{eq:stencil_sys} for each point of mesh. This can however be performed efficiently by noting that the $C, \bar{C}$ and $D$ matrices are derived solely from Finite Difference stencils on regular grids and thus do not change from point to point. These matrices are in fact identical across all meshes and may therefore be precomputed once and stored. The matrix $C$ may also be prefactored for instance into QR form to further streamline the linear system solve.

The only quantity then that needs to be computed for stencil assembly for each point of a given mesh is the matrix $X$, which can be cheaply assembled as it contains information only on nodal spatial coordinates of neighbouring points. This does not tend to incur any extra cost in memory as most mesh data structures already contain neighbour information for each point.

\subsection{Numerical differentiation example}
\label{sect:diff2d}
As an example, we consider a second order accurate approximation to the function $f(x,y) = e^{-x^2 -y^2}$ on a square, triangular and pentagonal meshes using 9-point regular/8-point extraordinary stencils shown in  Fig. \ref{fig:stencil_2d_centre3}. Symmetric central difference stencils are used on the interior, whereas one sided stencils are used on the boundaries to approximate the derivatives. Around extraordinary faces symmetric central 3 quadrant stencils are used. The initial planar, triangular, pentagonal meshes with zero refinements applied and the stencil points used for the example are shown also in Fig. \ref{fig:stencil_2d_centre3}.

We construct the function
\begin{equation}
    \begin{split}
         F_2(x,y;x_{k},y_{k}) &= f(x,y) - f(x_k, y_k) - (x-x_k)\frac{\partial f}{\partial x}(x_k,y_k) \\
         &- (y-y_k)\frac{\partial f}{\partial y}(x_k,y_k) - \frac{(x-x_k)^2}{2}\frac{\partial^2 f}{\partial x^2}(x_k,y_k) \\
         &- (x-x_k)(y-y_k)\frac{\partial^2 f}{\partial x \partial y}(x_k,y_k) \\
         &- \frac{(y-y_k)^2}{2}\frac{\partial^2 f}{\partial y^2}(x_k,y_k) \\
         \tilde{F}(\xi_1, \xi_2;x_{k},y_{k}) &= F_2\bigg( x(\xi_1,\xi_2), y(\xi_1,\xi_2); x_{k},y_{k} \bigg)
    \end{split}
\end{equation}
for which we compute the derivatives $\frac{\partial^j}{\partial \xi_1^{j_1} \partial \xi_2^{j_2}}, j_1 + j_2 = j \leq 2$ of $\tilde{F_p}$ at each vertex $(x_k, y_k)$ of a given mesh. Denoting the special case of a vertex at an extraordinary face with subscript $\text{EF}$, this leads to solving a linear system $CXDu = \bar{C} f$, which for interior non-boundary vertices each term is as follows:
\begin{equation*}
    \bar{C} =  \begin{bmatrix}
    0 & 0 & 0 & -1 & 0 & 1 & 0 & 0 & 0\\
    0 & -1 & 0 & 0 & 0 & 0 & 0 & 1 & 0\\
    0 & 0 & 0 & 1 & -2 & 1 & 0 & 0 & 0\\
    1 & 0 & -1 & 0 & 0 & 0 & -1 & 0 & 1\\
    0 & 1 & 0 & 0 & -2 & 0 & 0 & 1 & 0\\
    \end{bmatrix}
    ,
    \bar{C}_{\text{EF}} =  \begin{bmatrix}
            0 & 0 & -1 & 0 & 1 & 0 & 0 & 0\\
            -1 & 0 & 0 & 0 & 0 & 0 & 1 & 0\\
            0 & 0 & 1 & -2 & 1 & 0 & 0 & 0\\
            \frac{1}{2} & -\frac{1}{2} & \frac{1}{2} & -1 & \frac{1}{2} & -\frac{1}{2} & \frac{1}{2} & 0\\
            1 & 0 & 0 & -2 & 0 & 0 & 1 & 0\\
        \end{bmatrix}
\end{equation*}
are the full regular finite difference stencils with step size $h=1$ for the derivatives in order $\frac{\partial}{\partial x}, \frac{\partial}{\partial y}, \frac{\partial^2}{\partial x^2}, \frac{\partial^2}{\partial x \partial y}, \frac{\partial^2}{\partial y^2}$ for a regular point and a point on an extraordinary face respectively,
\begin{equation*}
    C =  \begin{bmatrix}
    0 & 0 & 0 & -1 & 1 & 0 & 0 & 0\\
    0 & -1 & 0 & 0 & 0 & 0 & 1 & 0\\
    0 & 0 & 0 & 1 & 1 & 0 & 0 & 0\\
    1 & 0 & -1 & 0 & 0 & -1 & 0 & 1\\
    0 & 1 & 0 & 0 & 0 & 0 & 1 & 0\\
    \end{bmatrix}
    ,
    C_{\text{EF}} =  \begin{bmatrix}
    0 & 0 & -1 & 1 & 0 & 0 & 0\\
    -1 & 0 & 0 & 0 & 0 & 1 & 0\\
    0 & 0 & 1 & 1 & 0 & 0 & 0\\
    \frac{1}{2} & -\frac{1}{2} & \frac{1}{2} & \frac{1}{2} & -\frac{1}{2} & \frac{1}{2} & 0\\
    1 & 0 & 0 & 0 & 0 & 1 & 0\\
    \end{bmatrix}
\end{equation*}
are the matrices constructed from $\bar{C}, \bar{C}_{\text{EF}}$ by deleting from each the column corresponding to the point $(x_k,y_k)$ at which the derivative is calculated,
\begin{equation*}
    D = \begin{bmatrix}
    \frac{1}{1!} & 0 & 0 & 0 & 0\\
    0 & \frac{1}{2!} & 0 & 0 & 0\\
    0 & 0 & \frac{1}{2!} & 0 & 0\\
    0 & 0 & 0 & \frac{1}{1!1!} & 0\\
    0 & 0 & 0 & 0 & \frac{1}{2!}\\
    \end{bmatrix}
\end{equation*}
a diagonal matrix of factorial denominators from Taylor Series' expansions,
\begin{equation*}
    X =  \begin{bmatrix}
    (x_1 - x_k) & (y_1 - y_k) & (x_1 - x_k)^2 & (x_1 - x_k)(y_1 - y_k) & (y_1 - y_k)^2\\
    \vdots & \ddots\\
    (x_{k-1} - x_k) & ... \\
    (x_{k+1} - x_k) & ... \\
    \vdots
    \end{bmatrix}
\end{equation*}
the Vandermonde matrix of the stencil point positions and vectors
\begin{equation*}
    f =
    \begin{bmatrix}
    f(x_{1}, y_{1})\\
    \vdots \\
    f(x_{9}, y_{9})
    \end{bmatrix}
    ,
    f_{\text{EF}} =
    \begin{bmatrix}
    f(x_{1}, y_{1})\\
    \vdots \\
    f(x_{8}, y_{8})
    \end{bmatrix}
\end{equation*}
the vectors of values of $f(x,y)$ at the stencil points $(x_1,y_1),...,(x_n,y_n)$.

\begin{figure}
    \centering
    \begin{tikzpicture}[font=\large, scale=0.8]
        \begin{semilogyaxis}[
            width=0.4\textwidth,
            height=0.5\textwidth,
            xlabel = {n},
            xlabel style={ yshift=-1ex },
            ylabel = {Max Error},
            ylabel style={ yshift=2ex },
            ymin=1e-6, ymax=1e-0,
            xmin=-0.2, xmax=4.2,
            xtick={0,1,2,3,4},
            xticklabels={0,1,2,3,4},
            ytick={1e0,1e-1,1e-2,1e-3,1e-4,1e-5,1e-6},
            yticklabels={$10^{0}$,$10^{-1}$,$10^{-2}$,$10^{-3}$,$10^{-4}$,$10^{-5}$, $10^{-6}$},
            grid = both,
            grid style = {line width=.1pt, draw=gray!15},
            major grid style = {line width=.2pt, draw=gray!50},
        ]
        \addplot[color=black]
        table[ x=n, y=dx ]{dat/diff2d_reg.csv};
        \label{plot:diff_d1}
            
        \addplot[color=black]
        table[ x=n, y=dy ]{dat/diff2d_reg.csv};
        \label{plot:diff_d2}
            
        \addplot[dashed]
        table[ x=n, y=dxx ]{dat/diff2d_reg.csv};
        \label{plot:diff_d3}
            
        \addplot[dashed]
        table[ x=n, y=dxy ]{dat/diff2d_reg.csv};
        \label{plot:diff_d4}
            
        \addplot[dashed]
        table[ x=n, y=dyy ]{dat/diff2d_reg.csv};
        \label{plot:diff_d4}
            
        \addplot[domain=2.5:3.8, samples=2] {0.015*3^(-2*x)};
        \node at (axis cs:3.1, 6e-6) {2};
            
        \end{semilogyaxis}
    \end{tikzpicture}%
    ~
    \begin{tikzpicture}[font=\large, scale=0.8]
        \begin{semilogyaxis}[
            width=0.4\textwidth,
            height=0.5\textwidth,
            xlabel = {n},
            xlabel style={ yshift=-1ex },
            ymin=1e-6, ymax=1e-0,
            xmin=-0.2, xmax=4.2,
            xtick={0,1,2,3,4},
            xticklabels={0,1,2,3,4},
            ytick={1e0,1e-1,1e-2,1e-3,1e-4,1e-5,1e-6},
            yticklabels={$10^{0}$,$10^{-1}$,$10^{-2}$,$10^{-3}$,$10^{-4}$,$10^{-5}$, $10^{-6}$},
            grid = both,
            grid style = {line width=.1pt, draw=gray!15},
            major grid style = {line width=.2pt, draw=gray!50},
        ]
        \addplot[color=black]
        table[ x=n, y=dx ]{dat/diff2d_tri.csv};
        \label{plot:diff_d1}
            
        \addplot[color=black]
        table[ x=n, y=dy ]{dat/diff2d_tri.csv};
        \label{plot:diff_d2}
            
        \addplot[dashed]
        table[ x=n, y=dxx ]{dat/diff2d_tri.csv};
        \label{plot:diff_d3}
            
        \addplot[dashed]
        table[ x=n, y=dxy ]{dat/diff2d_tri.csv};
        \label{plot:diff_d4}
            
        \addplot[dashed]
        table[ x=n, y=dyy ]{dat/diff2d_tri.csv};
        \label{plot:diff_d4}
            
        \addplot[domain=2.5:3.8, samples=2] {0.015*3^(-2*x)};
        \node at (axis cs:3.1, 6e-6) {2};
        
        \addplot[domain=2.5:3.8, samples=2] {0.28*3^(-x)};
        \node at (axis cs:3.1, 2.5e-2) {1};
        \end{semilogyaxis}
    \end{tikzpicture}%
    ~
    \begin{tikzpicture}[font=\large, scale=0.8]
        \begin{semilogyaxis}[
            width=0.4\textwidth,
            height=0.5\textwidth,
            xlabel = {n},
            xlabel style={ yshift=-1ex },
            ymin=1e-6, ymax=1e-0,
            xmin=-0.2, xmax=4.2,
            xtick={0,1,2,3,4},
            xticklabels={0,1,2,3,4},
            ytick={1e0,1e-1,1e-2,1e-3,1e-4,1e-5,1e-6},
            yticklabels={$10^{0}$,$10^{-1}$,$10^{-2}$,$10^{-3}$,$10^{-4}$,$10^{-5}$, $10^{-6}$},
            grid = both,
            grid style = {line width=.1pt, draw=gray!15},
            major grid style = {line width=.2pt, draw=gray!50},
        ]
        \addplot[color=black]
        table[ x=n, y=dx ]{dat/diff2d_pent.csv};
        \label{plot:diff_d1}
            
        \addplot[color=black]
        table[ x=n, y=dy ]{dat/diff2d_pent.csv};
        \label{plot:diff_d2}
            
        \addplot[dashed]
        table[ x=n, y=dxx ]{dat/diff2d_pent.csv};
        \label{plot:diff_d3}
            
        \addplot[dashed]
        table[ x=n, y=dxy ]{dat/diff2d_pent.csv};
        \label{plot:diff_d4}
            
        \addplot[dashed]
        table[ x=n, y=dyy ]{dat/diff2d_pent.csv};
        \label{plot:diff_d4}
            
        \addplot[domain=2.5:3.8, samples=2] {0.01*3^(-2*x)};
        \node at (axis cs:3.1, 4e-6) {2};
        
        \addplot[domain=2.5:3.8, samples=2] {0.2*3^(-x)};
        \node at (axis cs:3.1, 1.5e-2) {1};
        
        \end{semilogyaxis}
    \end{tikzpicture}
    \caption{Numerical differentiation of $f(x)=e^{-x^2 - y^2}$ on the three meshes shown in Fig. \ref{fig:stencil_2d_centre3} under refinement. $n$ denotes number of refinement steps applied to the mesh. From left to right: regular grid, triangular mesh, pentagonal mesh. The step size $h$ on the regular grid is chosen to be equal to the average edge length on the unstructured meshes at zero refinements. Solid lines show maximum absolute error in $\frac{\partial f}{\partial x}, \frac{\partial f}{\partial y}$, dashed lines show maximum absolute error in $\frac{\partial^2 f}{\partial x^2}, \frac{\partial^2 f}{\partial x \partial y}, \frac{\partial^2 f}{\partial y^2}$. }
    \label{fig:diff2d}
\end{figure}
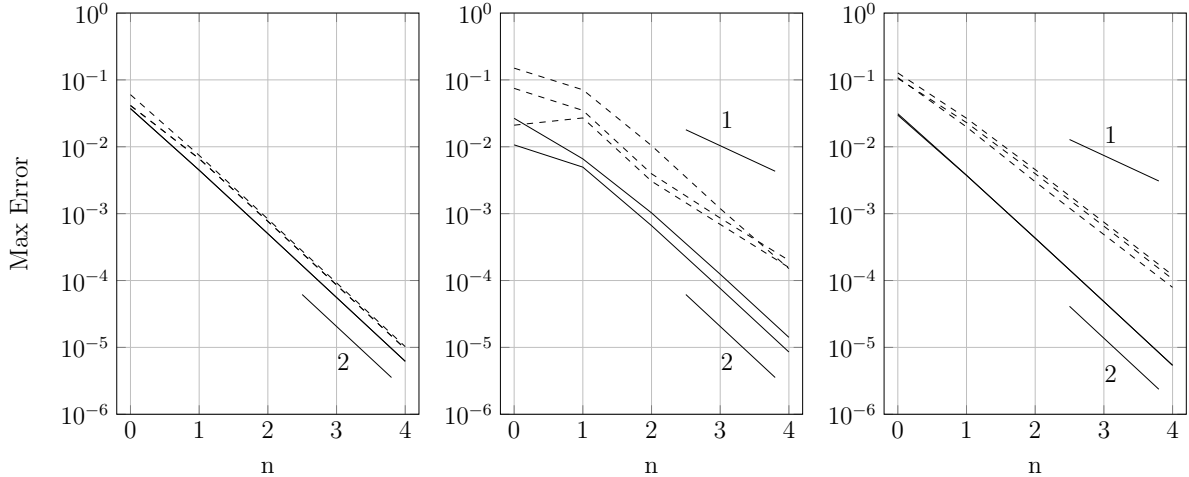

Solving the linear system for every non-boundary point of the meshes shown we obtain approximations to the derivatives $f(x,y) = e^{(-x^2 -y^2)}$. The maximum error of the approximations to each of the derivatives under refinement are shown in Fig. \ref{fig:diff2d}. In the case of a regular grid, the error in the derivatives converge with second order accuracy, whereas for the triangular and pentagonal meshes the error in the derivatives $\frac{\partial^2}{\partial x^2}, \frac{\partial^2}{\partial x \partial y}, \frac{\partial^2}{\partial y^2}$ converge at a slower rate although at a rate slightly higher than the predicted first order rate of convergence.

\section{Application to PDEs}
\begin{figure}
    \centering
    \includegraphics[scale=0.5]{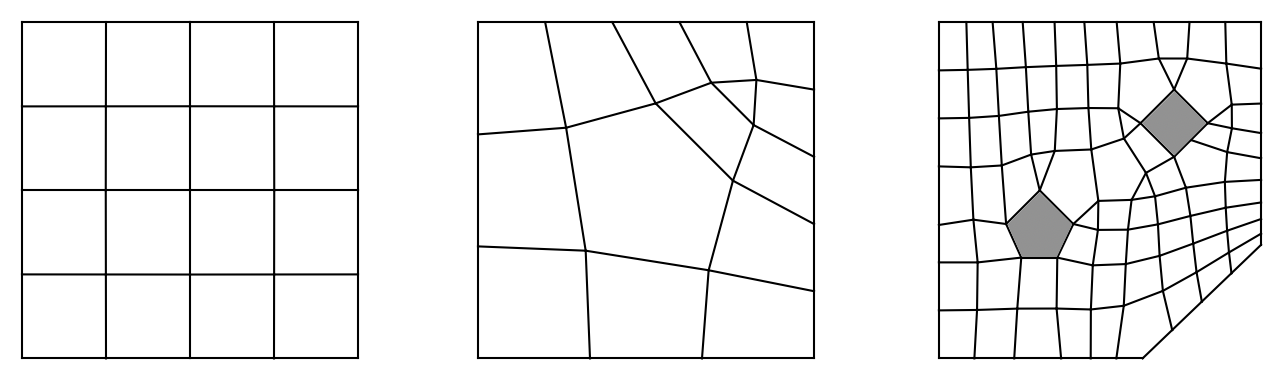}
    \caption{Meshes for PDE numerical examples at zero refinements. From left to right: regular structured grid plane mesh, unstructured plane mesh, unstructured polygonal mesh with two holes shown coloured in.}
    \label{fig:example_meshes}
\end{figure}

The numerical differencing procedure in the previous sections is applied the solution of various partial differential equations. For simplicity, in this paper we consider only examples in two-dimensions, although the theory is general and may be extended to higher dimensions in space. For the following examples, we consider the meshes at zero refinements shown in Fig. \ref{fig:example_meshes}: a plane discretised using a structured regular grid, a plane discretised using an unstructured mesh with two extraordinary faces, and a unstructured polygonal mesh with two holes containing multiple extraordinary faces.

\subsection{Poisson's equation in 2D}
\label{sect:poisson_ex}
We look to solve the boundary value problem in 2D
\begin{equation}
    \nabla^2 u = f
\end{equation}
on a domain $\Omega$ satisfying Dirichlet and Neumann boundary conditions
\begin{equation}
    u(x) = g(x) \text{ on } \partial \Omega_D ~,~ \frac{\partial u
    }{\partial n}(x) = h(x) \text{ on } \partial \Omega_N
\end{equation}
which in Cartesian coordinates can be expressed as
\begin{equation}
    \frac{\partial^2 u}{\partial x^2} + \frac{\partial^2 u}{\partial y^2} = f
\end{equation}
As a first example, we solve on the plane meshes discretised using a regular grid and an irregular grid with two extraordinary faces shown in Fig. \ref{fig:example_meshes}. The right hand side function is set as $f = 0$, and suitable Dirichlet boundary conditions set on the two horizontal boundaries and Neumann boundary conditions set on the two vertical boundaries such that the analytical solution $u(x) = \frac{1}{\sinh{\pi}}\sin{ (x \cdot \pi) } \sinh{ (y \cdot \pi) }$. Spatial derivatives are approximated using the 9-point regular point/8-point extraordinary face point stencils used in the example in Sect. \ref{sect:diff2d}, and Neumann boundary conditions implemented using the asymmetric boundary stencils shown in Fig. \ref{fig:stencil_2d_centre3}.

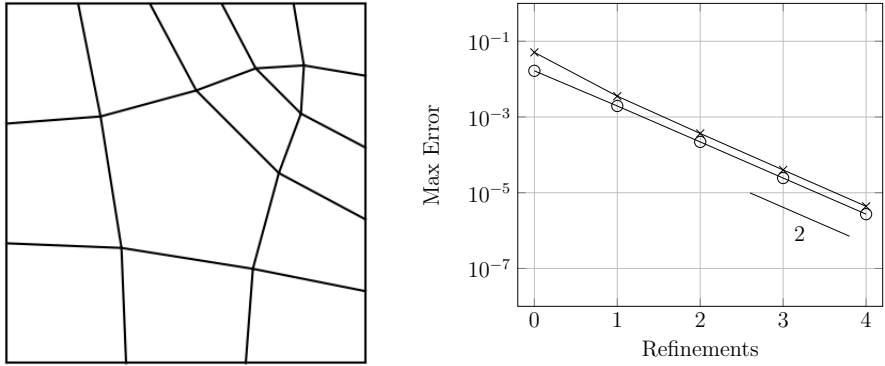
\begin{figure}[h]
    \centering
    \begin{minipage}{0.3\textwidth}
    \centering
    \includegraphics[width=\textwidth]{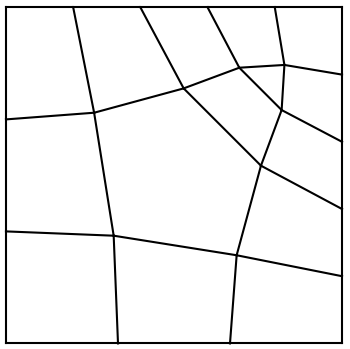}
    \end{minipage}
    \ \ \ \ %
    \begin{minipage}{0.38\textwidth}
    \resizebox{\linewidth}{!}{%
    \begin{tikzpicture}[font=\large]
        \begin{semilogyaxis}[
            xlabel = {Refinements},
            ylabel = {Max Error},
            ylabel style={ yshift=2ex },
            ymin=1e-8, ymax=1e-0,
            xmin=-0.2, xmax=4.2,
            xtick={0,1,2,3,4},
            xticklabels={0,1,2,3,4},
            ytick={1e-1,1e-3,1e-5,1e-7},
            yticklabels={$10^{-1}$,$10^{-3}$,$10^{-5}$,$10^{-7}$},
            grid = both,
            grid style = {line width=.1pt, draw=gray!15},
            major grid style = {line width=.2pt, draw=gray!50},
        ]
        \addplot[mark=x, color=black, mark size=3pt]
        table[ x=n, y=unstructured ]{dat/poi2d_3pt.csv};
        \label{plot:lw_unstructured}
        
        \addplot[mark=o, color=black, mark size=3pt]
        table[ x=n, y=structured ]{dat/poi2d_3pt.csv};
        \label{plot:lw_structured}
        
        \addplot[domain=2.6:3.8, samples=2] {0.003*3^(-2*x)};
        \node at (axis cs:3.2, 8e-7) {2};
        
        \end{semilogyaxis}
    \end{tikzpicture}
    }
    \end{minipage}
    \caption{Max error for Poisson's equation. \ref{plot:lw_structured} shows the error solving on structured grid plane mesh, \ref{plot:lw_unstructured} shows the error solving on unstructured plane mesh shown on the left at zero refinements.}
    \label{fig:poisson2d_3pt}
\end{figure}

Results for this example are shown in Fig. \ref{fig:poisson2d_3pt}. Denoting the computed solution as $\hat{u}(x)$, second order convergence in the max norm of the absolute error $|\hat{u}(x) - u(x)|$ is observed for both sets of meshes under refinement. The rate of convergence on the irregular mesh is higher than predicted by the local truncation errors, and higher than that observed in the numerical differentiation example in Sect. \ref{fig:diff2d}.

\begin{figure}[h]
    \centering
    \includegraphics[width=.95\textwidth]{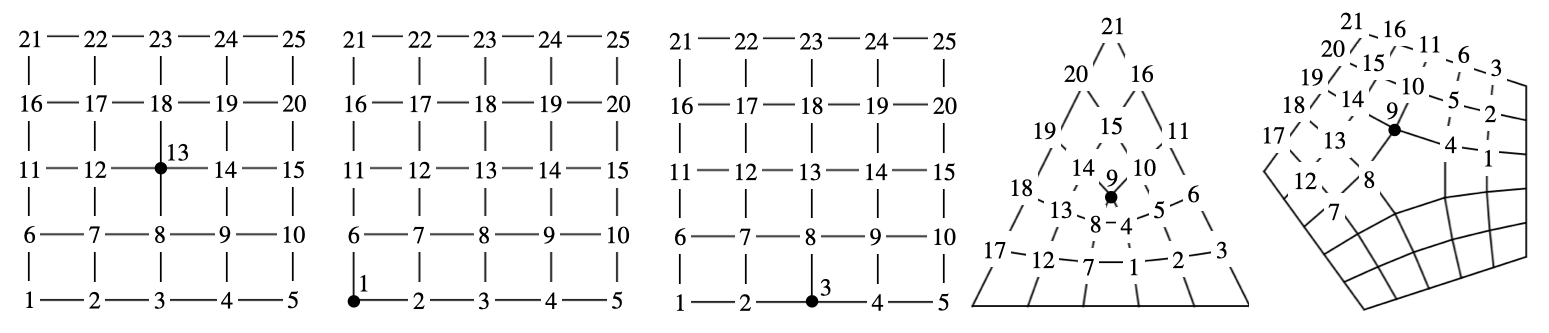}
    \caption{Stencil points and their numbering for Poisson's equation and biharmonic equation examples. Solid coloured vertex indicates point $(x_k,y_k)$ at which the function derivative is approximated. From left to right: symmetric stencil for regular point, one sided stencil for boundary point, one sided stencil for corner point, symmetric 3 quadrant stencil at an extraordinary face on a triangular mesh, symmetric 3 quadrant stencil at an extraordinary face on a pentagonal mesh.}
    \label{fig:stencil_2d_centre5}
\end{figure}

As a second example, we solve the boundary value problem on the unstructured two hole mesh shown in Sect. \ref{fig:diff2d}. As with the previous example the right hand side is set to be $f=0$, and boundary conditions set such that the analytical solution $u(x) = \frac{1}{\sinh{\pi}}\sin{ (x \cdot \pi) } \sinh{ (y \cdot \pi) }$. 

To obtain a higher order accurate solution we utilise the expanded 25-point regular point/21-point extraordinary face point stencils shown in Fig. \ref{fig:stencil_2d_centre5} to discretise the spatial derivatives of the PDE. Neumann boundary conditions are set on the two vertical boundaries which are implemented using the asymmetric stencils shown in Fig. \ref{fig:stencil_2d_centre5}, and Dirichlet boundary conditions set on all other boundaries. As a consequence of the choice of stencil, a requirement that extraordinary faces be separated by at least two rings regular four-sided faces must be satisfied for the ability to choose this stencil globally on the mesh. This means that the two holed mesh with zero refinement steps applied shown in Fig. \ref{fig:example_meshes} must be refined at least once for this stencil to be applied in this example.

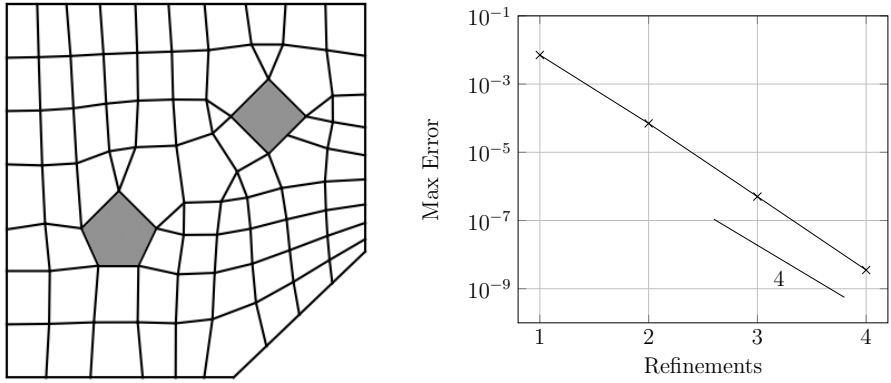
\begin{figure}[h]
    \centering
    \begin{minipage}{0.3\textwidth}
    \centering
    \includegraphics[width=\textwidth]{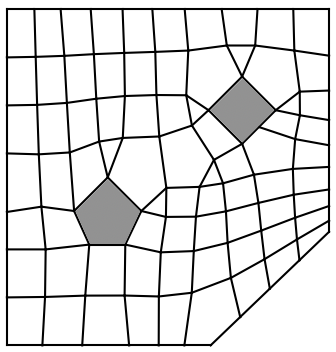}
    \end{minipage}
    \ \ \ \ %
    \begin{minipage}{0.38\textwidth}
    \resizebox{\linewidth}{!}{%
    \begin{tikzpicture}[font=\large]
        \begin{semilogyaxis}[
            xlabel = {Refinements},
            ylabel = {Max Error},
            ylabel style={ yshift=2ex },
            ymin=1e-10, ymax=1e-1,
            xmin=0.8, xmax=4.2,
            xtick={1,2,3,4},
            xticklabels={1,2,3,4},
            grid = both,
            grid style = {line width=.1pt, draw=gray!15},
            major grid style = {line width=.2pt, draw=gray!50},
        ]
        \addplot[mark=x, color=black, mark size=3pt]
        table[ x=n, y=unstructured ]{dat/poi2d_twohole.csv};
        
        \addplot[domain=2.6:3.8, samples=2] {0.01*3^(-4*x)};
        \node at (axis cs:3.2, 2e-9) {4};

        \end{semilogyaxis}
    \end{tikzpicture}%
    }
    \end{minipage}
    \caption{Max error for Poisson's equation on unstructured polygonal mesh with two holes shown on the left at zero refinements.}
    \label{fig:poisson2d_5pt}
\end{figure}

Results for this example are shown in Fig. \ref{fig:poisson2d_5pt}. Denoting the exact solution as $\hat{u}(x)$, we observe error $|\hat{u}(x) - u(x)|$ in the max norm to be fourth order accurate. Despite the presence of numerous extraordinary faces the order of accuracy observed is equal to that of a regular grid discretisation rather than being one order lower as suggested by the local truncation error.

\subsection{Biharmonic equation in 2D}
We apply the method to solve the fourth order biharmonic equation
\begin{equation}
    \nabla^4 u = f
\end{equation}
on a domain $\Omega$ satisfying boundary conditions
\begin{equation}
    u(x) = g(x) \text{ on } \partial \Omega ~,~ \frac{\partial u
    }{\partial n}(x) = h(x) \text{ on } \partial \Omega
\end{equation}
which can be written in Cartesian coordinates as
\begin{equation}
    \frac{\partial^4 u}{\partial x^4} + \frac{\partial^4 u}{\partial x^2 \partial y^2} + \frac{\partial^4 u}{\partial y^4} = f
\end{equation}

We compute a solution to the equation for the two planar domains discretised using a regular grid and an irregular mesh shown in Fig. \ref{fig:example_meshes}. Spatial derivatives are computed using the 25-point regular point/21-point extraordinary face point stencils shown in Fig. \ref{fig:stencil_2d_centre5}. Similar to the previous example with Poisson's equation, for this stencil to be used globally each extraordinary face must be separated by at least two rings of regular four-sided faces, meaning on refinement step must be applied to the meshes shown before the stencil is applied.

For this example, the right hand side is set to be $f = 0$. Both Dirichlet and Neumann boundary conditions are set on all boundaries such that the solution $u(x,y) = \frac{1}{\sinh{(\pi)}} (x^2 + y^2) \sin{(\pi x)} \sinh{(\pi y)}$.

\begin{figure}[h]
    \centering
    \begin{minipage}{0.3\textwidth}
    \centering
    \includegraphics[width=\textwidth]{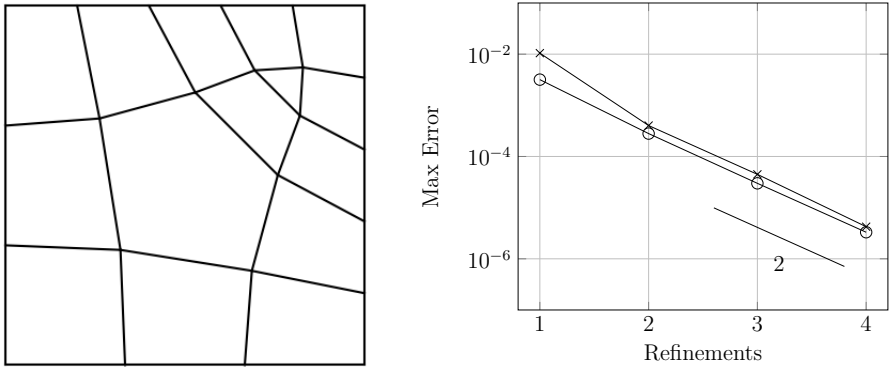}
    \end{minipage}
    \ \ \ \ %
    \begin{minipage}{0.38\textwidth}
    \resizebox{\linewidth}{!}{%
    \begin{tikzpicture}[font=\large]
        \begin{semilogyaxis}[
            xlabel = {Refinements},
            ylabel = {Max Error},
            ylabel style={ yshift=2ex },
            ymin=1e-7, ymax=1e-1,
            xmin=0.8, xmax=4.2,
            xtick={1,2,3,4},
            xticklabels={1,2,3,4},
            ytick={1e-2,1e-4,1e-6},
            yticklabels={$10^{-2}$,$10^{-4}$,$10^{-6}$},
            grid = both,
            grid style = {line width=.1pt, draw=gray!15},
            major grid style = {line width=.2pt, draw=gray!50},
        ]
        \addplot[mark=x, color=black, mark size=3pt]
        table[ x=n, y=unstructured ]{dat/biharmonic2d_5pt.csv};
        
        \addplot[mark=o, color=black, mark size=3pt]
        table[ x=n, y=structured ]{dat/biharmonic2d_5pt.csv};
        
        \addplot[domain=2.6:3.8, samples=2] {0.003*3^(-2*x)};
        \node at (axis cs:3.2, 8e-7) {2};
        
        \end{semilogyaxis}
    \end{tikzpicture}%
    }
    \end{minipage}
    \caption{Max error for biharmonic equation. \ref{plot:lw_structured} shows the error solving on structured plane mesh, \ref{plot:lw_unstructured} shows the error solving on unstructured plane mesh shown on the left at zero refinements.}
    \label{fig:biharmonic2d}
\end{figure}

Results for this example are shown in Fig. \ref{fig:biharmonic2d}. Denoting the computed solution as $\hat{u}(x)$, we observe the absolute error $|\hat{u}(x) - u(x)|$ in the max norm to be second order accurate for both sets of meshes. This observation is consistent with that from the previous example with Poisson's equation, where the order of convergence on an unstructured mesh is an order higher than suggested by the local truncation error and equal to that of a regular grid discretisation.

\subsection{Minimal Surfaces}
\begin{figure}[h]
    \centering
    \begin{minipage}{0.3\textwidth}
    \includegraphics[width=\textwidth]{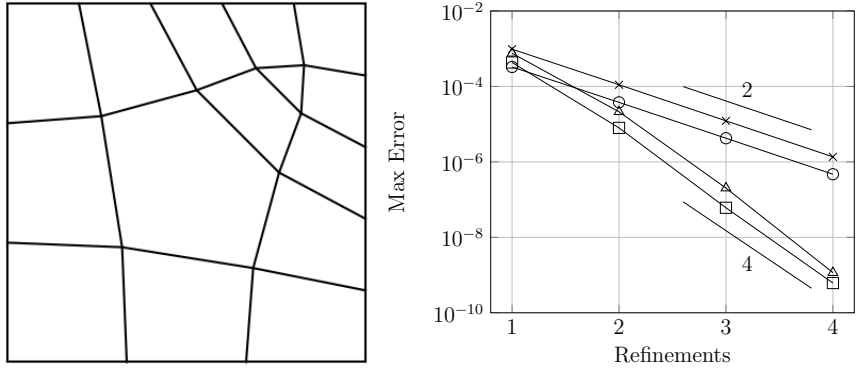}
    \end{minipage}
    \begin{minipage}{0.38\textwidth}
    \resizebox{\linewidth}{!}{%
    \begin{tikzpicture}[font=\large]
        \begin{semilogyaxis}[
            xlabel = {Refinements},
            ylabel = {Max Error},
            ylabel style={ yshift=2ex },
            ymin=1e-10, ymax=1e-2,
            xmin=0.8, xmax=4.2,
            xtick={1,2,3,4},
            xticklabels={1,2,3,4},
            ytick={1e-2,1e-4,1e-6,1e-8,1e-10},
            yticklabels={$10^{-2}$,$10^{-4}$,$10^{-6}$,$10^{-8}$,$10^{-10}$},
            grid = both,
            grid style = {line width=.1pt, draw=gray!15},
            major grid style = {line width=.2pt, draw=gray!50},
        ]
        \addplot[mark=x, color=black, mark size=3pt]
        table[ x=n, y=unstructured2 ]{dat/minimal_surface.csv};
        \label{plot:minsurf_u2}
        
        \addplot[mark=o, color=black, mark size=3pt]
        table[ x=n, y=structured2 ]{dat/minimal_surface.csv};
        \label{plot:minsurf_s2}
        
        \addplot[mark=triangle, color=black, mark size=3pt]
        table[ x=n, y=unstructured4 ]{dat/minimal_surface.csv};
        \label{plot:minsurf_u4}
        
        \addplot[mark=square, color=black, mark size=3pt]
        table[ x=n, y=structured4 ]{dat/minimal_surface.csv};
        \label{plot:minsurf_s4}
        
        \addplot[domain=2.6:3.8, samples=2] {0.03*3^(-2*x)};
        \node at (axis cs:3.2, 8e-5) {2};
        
        \addplot[domain=2.6:3.8, samples=2] {0.008*3^(-4*x)};
        \node at (axis cs:3.2, 2e-9) {4};
        
        \end{semilogyaxis}
    \end{tikzpicture}%
    }
    \end{minipage}
    \caption{Max error for minimal surface example. The unstructured planar mesh considered is shown on the left at zero refinements. The errors using the 9-point/8-point extraordinary face stencil are shown by \ref{plot:minsurf_s2} for the structured grid plane mesh, and \ref{plot:minsurf_u2} for the unstructured plane mesh. The errors using the 25-point/21-point extraordinary face stencil are shown by \ref{plot:minsurf_s4} for the structured grid plane mesh, and \ref{plot:minsurf_u4} for the unstructured plane mesh.}
    \label{fig:minimal_surface}
\end{figure}

We consider the problem of finding a minimal surface, defined to be a surface that satisfies the constraint zero mean curvature globally. Minimal surfaces arise in numerous physical applications such as in soap films, and are governed by the nonlinear equation
\begin{equation}
    \bigg( 1+\frac{\partial u}{\partial x}^2 \bigg) \frac{\partial^2 u}{\partial y^2} - 2 \frac{\partial u}{\partial x} \frac{\partial u}{\partial y} \frac{\partial^2 u}{\partial x \partial y} + \bigg( 1+\frac{\partial u}{\partial y}^2 \bigg) \frac{\partial^2 u}{\partial x^2} = 0
\end{equation}
Dirichlet boundary conditions are set for this example such that the solution is the Scherk surface, which can be expressed in the form $u(x,y) = c \bigg( \ln{ \big( \cos{(\frac{y}{c})} \big) } - \ln{ \big( \cos{(\frac{x}{c})} \big) } \bigg)$ where $c$ is some positive real number. We choose $c=1$ and solve the problem on a planar domain $\Omega = [-1,1]^2$ discretised using the regular and irregular planar meshes shown in Fig. \ref{fig:example_meshes}. For the choice of stencil, we look at both the 9-point regular point/8-point extraordinary face point stencil and the 25-point regular point/21-point extraordinary face point stencil shown in Sect. \ref{sect:diff2d} and Sect. \ref{sect:poisson_ex} respectively, the results of which are shown in Fig. \ref{fig:minimal_surface}.

Denoting the computed solution as before as $\hat{u}(x)$, we observe the error in the max norm $|\hat{u}(x) - u(x)|$ for both meshes using the 9-point/8-point extraordinary face stencil to be second order accurate, while the error in the max norm using the 25-point/21-point extraordinary face stencil is fourth order accurate. As with the previous examples we observe the order of convergence to be unaffected by the presence of topological defects or other irregularities in the mesh.

\subsection{Advection equation with upwinding}
We consider the time-dependent scalar advection equation
\begin{equation}
   \frac{\partial u}{\partial t} + v \cdot \nabla u = 0 , ~~u(x,0) = u_0(x)
\end{equation}
where $v$ is the advection velocity. It is well known that on a regular grid, use of a centred stencil results in a skew-symmetric matrix with purely imaginary eigenvalues, which restricts the range of timestepping methods that may be employed to solve the system. To stabilise the system such that the system eigenvalues all have negative real part, various discretisation strategies have been introduced including the Lax-Wendroff and Beam-Warming methods to name a few. We focus here on extending upwinding methods to unstructured meshes for solving the scalar advection equation.

\begin{figure}
    \centering
    \includegraphics[scale=0.45]{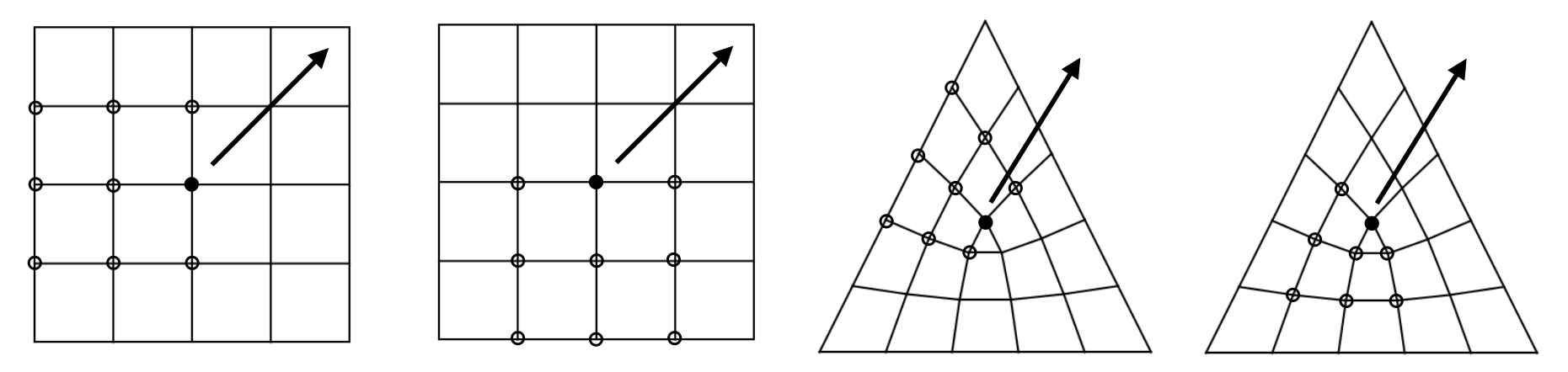}
    \caption{Stencils for x and y gradient components in 2D. Black arrow shows direction of velocity at bolded point. Left two meshes show points used to construct stencils for $x$ and $y$ derivatives respectively. Right two meshes show points used to construct stencils for $x$ and $y$ at an extraordinary face.}
    \label{fig:upwind_stencil}
\end{figure}

Upwinding methods are inspired by the method of characteristics for analysing hyperbolic PDEs, wherein the stencil is chosen to be skewed downwind from the direction of the advection velocity. In 1D, this simply means that for a positive velocity, a backwards difference stencil such as Equation \ref{eq:1d_backward} is chosen for approximating the first derivatives of $u$. For a regular grid in higher dimensions this idea can be extended for each of the directions of the velocity independently as shown in Fig. \ref{fig:upwind_stencil}.

To extend this idea to unstructured meshes, at each point of the mesh for each direction of the velocity we identify the neighbouring point closest to being downwind to the component of the velocity. A stencil is then constructed such that it is centred on this vertex for that component of the gradient. An extra complication is introduced at topological defects at extraordinary faces, where a regular grid stencil cannot be defined uniquely. In this case an analogous procedure to the one done above can be performed wherein one quadrant is omitted to form the stencil. For this example 9-point regular/8-point irregular stencils are chosen such that the method is second order accurate. The points used to construct these stencils in two dimensions are shown in Fig. \ref{fig:upwind_stencil}.

\begin{figure}[h]
    \centering
    \begin{minipage}{0.38\textwidth}
    \resizebox{\linewidth}{!}{%
    \begin{tikzpicture}[font=\large]
        \begin{axis}[
            xlabel = {Re},
            ylabel = {Im},
            ylabel style={ yshift=2ex },
            ymin=-70, ymax=70,
            xmin=-200, xmax=10,
            grid = both,
            grid style = {line width=.1pt, draw=gray!15},
            major grid style = {line width=.2pt, draw=gray!50},
        ]
        \addplot[mark=x, only marks, color=black, mark size=1pt]
        table[ x=unstruct_real, y=unstruct_imag ]{dat/upwindconvect_eig.csv};
        \label{plot:upwindconvect_unstruct_eig}
        
        \end{axis}
    \end{tikzpicture}%
    }
    \end{minipage}%
    \begin{minipage}{0.343\textwidth}
    \resizebox{\linewidth}{!}{%
    \begin{tikzpicture}[font=\large]
        \begin{axis}[
            xlabel = {Re},
            ylabel style={ yshift=2ex },
            ymin=-70, ymax=70,
            xmin=-100, xmax=10,
            grid = both,
            grid style = {line width=.1pt, draw=gray!15},
            major grid style = {line width=.2pt, draw=gray!50},
        ]
        \addplot[mark=o, only marks, color=black, mark size=1pt]
        table[ x=struct_real, y=struct_imag ]{dat/upwindconvect_eig.csv};
        \label{plot:upwindconvect_struct_eig}
        
        \end{axis}
    \end{tikzpicture}%
    }
    \end{minipage}
    \caption{Eigenvalues of upwinded advection matrix. Left shows eigenvalues of unstructured plane mesh, right shows eigenvalues of structured grid, each under 2 refinements.}
    \label{fig:upwind_convect_eig}
\end{figure}
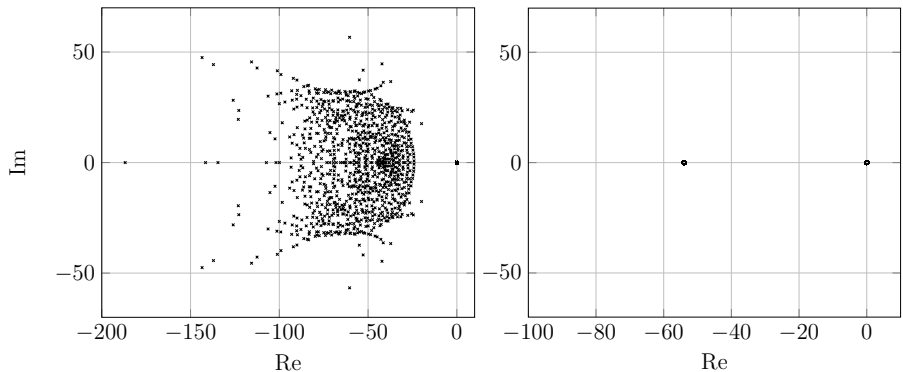

A square domain $[-1,1]^2$ is chosen with initial condition $u(x,y) = \exp[-20((x+0.3)^2 + (y+0.3)^2 )]$ and velocity vector $v(x,y) = (1,1)$. Time integration was performed using implicit Euler for simplicity with timestep $10^{-4}$ and the simulation run until $T = 1.0$. The domain is discretised using both a regular grid and an irregular mesh. To check for stability of the method, the eigenvalues of the system matrix under 2 refinements is shown in Fig. \ref{fig:upwind_convect_eig}.

The eigenvalues for the unstructured system while all have negative real part fluctuate in magnitude compared to those on the regular grid. This reflects the fact that distances between points on the unstructured mesh vary, especially around extraordinary faces where edge lengths shrink at a different rate under refinement compared to at quadriliaterals. This is a well known phenomenon in the field of subdivision surfaces known as the characteristic map. 

\begin{figure}[h]
    \centering
    \begin{minipage}{0.3\textwidth}
    \includegraphics[width=\textwidth]{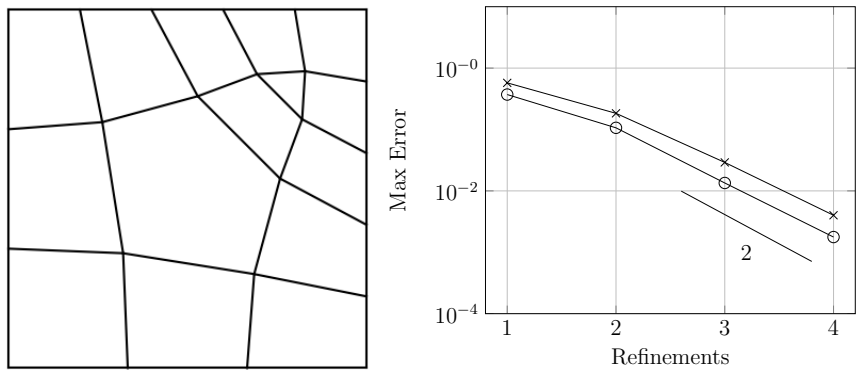}
    \end{minipage}
    \begin{minipage}{0.38\textwidth}
    \resizebox{\linewidth}{!}{%
    \begin{tikzpicture}[font=\large]
        \begin{semilogyaxis}[
            xlabel = {Refinements},
            ylabel = {Max Error},
            ylabel style={ yshift=2ex },
            ymin=1e-4, ymax=1e1,
            xmin=0.8, xmax=4.2,
            xtick={1,2,3,4},
            xticklabels={1,2,3,4},
            ytick={1e-0,1e-2,1e-4},
            yticklabels={$10^{-0}$,$10^{-2}$,$10^{-4}$},
            grid = both,
            grid style = {line width=.1pt, draw=gray!15},
            major grid style = {line width=.2pt, draw=gray!50},
        ]
        \addplot[mark=x, color=black, mark size=3pt]
        table[ x=nRefine, y=unstruct ]{dat/upwindconvect.csv};
        \label{plot:upwindconvect_unstruct}
        
        \addplot[mark=o, color=black, mark size=3pt]
        table[ x=nRefine, y=struct ]{dat/upwindconvect.csv};
        \label{plot:upwindconvect_struct}
        
        \addplot[domain=2.6:3.8, samples=2] {3*3^(-2*x)};
        \node at (axis cs:3.2, 1e-3) {2};
        
        \end{semilogyaxis}
    \end{tikzpicture}%
    }
    \end{minipage}
    \caption{Max error for upwinded convection. \ref{plot:upwindconvect_struct} shows the error on structured plane mesh, \ref{plot:upwindconvect_unstruct} shows the error on unstructured plane mesh.}
    \label{fig:upwind_convect}
\end{figure}

Results for this example are shown in Fig. \ref{fig:upwind_convect}. Denoting the computed solution as $\hat{u}$, second order convergence of the error in the max norm $|\hat{u}(x) - u(x)|$ for both sets of meshes is observed. This suggests that the unstructured upwinding procedure retains second order accuracy despite the presence of extraordinary faces in the mesh.

\section{Conclusion}

We have introduced a new framework for deriving finite difference discretisations on irregular grids. By defining an auxiliary function with high degrees of regularity, we used standard stencils for equally-spaced points to find the high-order equivalents on the original non-uniform grid. We extended the method to higher spatial dimensions on fully unstructured meshes of quadrilateral elements, using a subdivision-based refinement strategy and defining the node-points on the dual mesh. We used these results to demonstrate the high-order accuracy of the method for various PDEs and meshes.

In our future work, we will study other equations such as the Navier-Stokes equations, as well as how to incorporate state-of-the-art finite difference techniques using this formulation. We are also interested in 3D problems, which we believe will be a straight-forward extension (although generating the required fully unstructured hexahedral meshes is a well-known challenge). Finally, we will study the computational performance of the method, in particular when compared with corresponding high-order finite element / discontinuous Galerkin discretisations.

\section*{Acknowledgments}

This work was supported in part by the Director, Office of Science, Office of
Advanced Scientific Computing Research, U.S. Department of Energy under
Contract No. DE-AC02-05CH11231.

\bibliographystyle{plain}
\bibliography{references}

\end{document}